\documentclass[a4paper]{amsart}
\usepackage[utf8]{inputenc}
\usepackage[english]{babel}

\usepackage{amscd,amssymb}

\usepackage[bbgreekl]{mathbbol}
\DeclareSymbolFont{AMSb}{U}{msb}{m}{n}
\DeclareSymbolFontAlphabet{\mathbb}{AMSb}

\usepackage{amsthm}
\usepackage{graphicx}
\usepackage{fancybox}
\usepackage{epic,eepic}
\usepackage{amstext}
\usepackage{mathtools}
\usepackage{esint}
\usepackage[square,numbers]{natbib}
\usepackage{booktabs}
\usepackage[hide links]{hyperref}
\usepackage{comment}
\usepackage{mathtools}
\usepackage{tikz,pgfplots}
\usepackage{subcaption}

\usepackage[noabbrev, capitalise, nameinlink]{cleveref}
\crefname{equation}{}{}
\crefname{assumption}{Assumption}{Assumptions}
\crefformat{equation}{\textup{#2(#1)#3}}

\newtheorem{theorem}{Theorem}[section]

\newtheorem{lemma}[theorem]{Lemma}
\newtheorem{example}[theorem]{Example}
\theoremstyle{definition}

\newtheorem{assumption}[theorem]{Assumption}
\theoremstyle{remark}
\newtheorem{remark}[theorem]{Remark}
\numberwithin{theorem}{section}
\numberwithin{equation}{section}
\numberwithin{figure}{section}

\def\TH{\mathcal T_H}
\def\Th{\mathcal T_h}

\def\diam{\operatorname{diam}}

\def\div{\nabla\cdot}
\def\with{\,:\,}
\def\dx{\,\mathrm{d}x}
\def\ds{\,\mathrm{d}\sigma}
\def\tint{\begingroup\textstyle \int\endgroup}
\def\Nb{\mathsf{N}}
\def\blambda{\boldsymbol{\lambda}}
\def\bmu{\boldsymbol{\mu}}

\numberwithin{equation}{section}
\numberwithin{theorem}{section}

\AtBeginDocument{%
	\def\MR#1{}
}

\title[A High-Order LOD Method for Heterogeneous Stokes Problems]{A High-Order Localized Orthogonal Decomposition Method for Heterogeneous Stokes Problems}
\author[M.~Hauck, A. Lozinski]{Moritz Hauck$^*$, Alexei Lozinski$^\dagger$}
\address{${}^*$ Institute for Applied and Numerical Mathematics, Karlsruhe Institute of Technology, Englerstr.~2, 76131 Karlsruhe, Germany}
\email{moritz.hauck@kit.edu}
\address{${}^{\dagger}$ Université Marie et Louis Pasteur, CNRS, LmB (UMR 6623), F-25000 Besançon, France.}
\email{alexei.lozinski@univ-fcomte.fr}

\begin{document}

\begin{abstract} 
In this paper, we propose a high-order extension of the 
multiscale method introduced by the authors in [SIAM J. Numer. Anal., 63(4) (2025), pp. 1617–1641] for heterogeneous Stokes problems, while also providing several other improvements, including a better localization strategy and a more precise pressure reconstruction. The proposed method is based on the Localized Orthogonal Decomposition methodology and achieves optimal convergence orders under minimal structural assumptions on the coefficients. A key feature of our approach is the careful design of so-called quantities of interest, defining functionals of the solution whose values the multiscale approximation aims to reproduce exactly. Their selection is particularly delicate in the context of Stokes problems due to potential conflicts arising from the divergence-free constraint. 
We prove the exponential decay of the problem-adapted basis functions, justifying their localized computation in practical implementations. A rigorous a priori error analysis proves high-order convergence for both velocity and pressure, if the basis supports grow logarithmically with the desired accuracy. Numerical experiments confirm the theoretical findings. 
\end{abstract}



\maketitle

\section{Introduction}
We consider a heterogeneous Stokes problem posed on a bounded Lipschitz polytope $\Omega \subset \mathbb{R}^n$, $n \in \{2,3\}$. For a given external force $f$, the problem is to find a velocity $u$ and a pressure $p$ satisfying:
\begin{equation}
	\label{pbStokes}
	\left\{
	\begin{aligned}
		- \div (\nu \nabla u) + \sigma
		u + \nabla p = f, & \qquad \text{in } \Omega,\\
		\div u = 0, & \qquad \text{in } \Omega,\\
		u = 0, & \qquad \text{on } \partial \Omega,
	\end{aligned}\right.
\end{equation}
where $\nu$ and $\sigma$ denote the viscosity and damping coefficients, respectively. These coefficients encode the heterogeneity of the medium and may exhibit roughness or oscillations across multiple, possibly non-separated, length scales. Heterogeneous Stokes problems such as \cref{pbStokes} arise naturally in a variety of applications. In magma modeling, for example, the viscosity depends on temperature and may vary significantly across the domain, cf.~\cite{Gutirrez2010}. Another typical scenario occurs in slow flows around many small obstacles, cf.~\cite{angot99}. In this case, we set $\nu$ to the physical viscosity and $\sigma = 0$ in the fluid region, while inside the obstacles both coefficients take large values, effectively modeling solid inclusions.
	
The numerical approximation of heterogeneous problems such as \cref{pbStokes} by standard finite element methods (FEMs) typically suffers from reduced convergence rates and pronounced pre-asymptotic effects when the computational mesh does not resolve the fine-scale variations of the coefficients. Since globally resolving all microscopic details is computationally prohibitive, it is desirable to design numerical methods that yield accurate approximations even on coarse meshes that do not necessarily resolve the coefficients’ heterogeneities. This is realized through the construction of problem-adapted basis functions, an approach underlying many modern multiscale methods. For elliptic diffusion-type problems, examples include the Heterogeneous Multiscale Method~\cite{EE03,EE05,AbdEEV12}, the (Generalized) Multiscale FEM~\cite{BabO83,BabCO94,HowW97,BabL11,EfeGH13}, the Multiscale Spectral Generalized FEM~\cite{BabL11,MaSD22}, rough polyharmonic splines~\cite{OwhZB14}, the Localized Orthogonal Decomposition (LOD)~\cite{MalP14,HenP13}, and gamblets~\cite{Owh17}. More recently, refined localization strategies within the LOD framework have been proposed; see, e.g., \cite{HauP23,FreHKP24}. Comprehensive overviews of multiscale methods are available in the textbooks~\cite{OwhS19,MalP20} and the review article~\cite{AltHP21}.
	
Several of the above-mentioned multiscale methods for diffusion-type problems have been successfully adapted to Stokes problems, which present additional difficulties due to the divergence-free constraint. For slowly varying perforated media, we refer, for instance, to~\cite{Brown2013,Brown12013}. The Multiscale FEM based on Crouzeix--Raviart elements, originally introduced in~\cite{LeBris2014}, has been applied to Stokes flows in perforated domains~\cite{Muljadi2015,Jankowiak2024,Feng22,balaziPhD}. Also a variant of the Generalized Multiscale FEM for Stokes problems in perforated media was developed in~\cite{Chung2021}. More recently, a lowest-order multiscale method for the Stokes problem within the LOD framework was proposed in~\cite{Hauck2025}. Finally, we mention the generalization of the operator-adapted wavelet approach from~\cite{budninskiy2019operator}, which treats general problems with differential constraints, including divergence-freeness.
	
While the techniques discussed above typically exhibit first-order convergence, high-order multiscale methods have also been developed. For diffusion-type problems, such extensions have been proposed, for instance, for the Heterogeneous Multiscale Method~\cite{LiPT12,AbdB12} and for the Multiscale FEM \cite{AllB05,HesZZ14}. Hybrid multiscale methods, which reduce global degrees of freedom to element boundaries, have also gained popularity as a means to achieve high-order convergence; see, e.g.,~\cite{HarPV13,AraHPV13,CicEL19}. For Stokes problems, a high-order variant of the Multiscale FEM has been proposed in~\cite{Feng22,balaziPhD}.
However, all of these high-order approaches require certain smoothness assumptions on the domain, the coefficients, and/or the exact solution to achieve convergence rates beyond first order. In the presence of rough coefficients, as frequently encountered in applications, such as composite materials with abrupt transitions between material properties, these conditions are typically not satisfied, and the coefficients are often only in~$L^\infty$. Obtaining high-order convergence in this setting is nontrivial and requires the careful design of problem-tailored approximation spaces, together with appropriate orthogonality properties.
For diffusion-type problems, this is achieved in~\cite{Maier2021,Dong2023,Hauck2025b}, where high-order multiscale method based on ideas of the LOD and gamblets were developed.  In these works, quantities of interest (QOIs), which are functionals of the solution that the multiscale approximation aims to preserve exactly, were defined as integrals with piecewise polynomials spaces. For Stokes problems, applying this approach directly would lead to an ill-posed numerical method, since the divergence-free constraint is not accounted for, and also the construction of a high-order method is more involved.

The construction of high-order multiscale methods for heterogeneous Stokes problems based on the LOD methodology is the focus of the present article. To this end, we employ a reformulation of Stokes problem as in~\cite{Hauck2025}, where the velocity belongs to $(H^1_0(\Omega))^n$, with divergence piecewise constant on the underlying coarse mesh, and the pressure Lagrange multiplier is chosen to be piecewise constant on the same mesh. We then apply the LOD methodology to this reformulated problem. For the method of degree $m$, our choice of QOIs is based on a suitable decomposition of the space of polynomials of degree $m$ into those that can be expressed as gradients of scalar polynomials of degree $m+1$ and a complementary space. For a suitably chosen complement space, the QOIs consist of weighted element integrals, where the weights are given by its basis functions, together with weighted normal face integrals, using weights given by polynomials of degree up to~$m$ on the faces.
 This construction enables high-order convergence rates to be extracted from the source term through appropriate orthogonality relations. The resulting problem-adapted basis functions decay exponentially, which justifies their localized computation in practical implementations.
\medskip

The novel aspects of this article are summarized as follows:
\begin{itemize}
    \item We provide a rigorous a priori error analysis of the proposed method, establishing its high-order convergence under minimal assumptions on the coefficients.The velocity approximation converges at order $m+2$ in the $H^1$-norm and at order $m+3$ in the $L^2$-norm, provided the number of element layers in the patches to which the basis functions are localized grows logarithmically with the desired accuracy. Moreover, the piecewise constant pressure approximation is shown to converge exponentially to the corresponding pressure averages as the number of element layers is increased.
    \item A novel post-processing step is introduced that reconstructs a pressure approximation of order $m+2$ in the $L^2$-norm. Compared to the post-processing in \cite{Hauck2025}, this approach is more sophisticated, incorporating appropriate coarse-scale piecewise polynomial corrections. In the lowest-order case $m=0$, it achieves second-order convergence, compared to first-order convergence for the lowest-order method in~\cite{Hauck2025}.
    \item Special attention in the analysis of the method is paid to the fine-scale discretization, which is essential for the practical computation of the local, infinite-dimensional problems defining its basis functions. Taylor--Hood and Scott--Vogelius finite elements are employed for this purpose. We prove the well-posedness of the fully discrete problems for the basis functions and that the convergence results remain valid after fine-scale discretization.
    \item The stabilized localization strategy introduced in \cite{Hauck2025b} for diffusion-type problems is adopted to the current setting of Stokes problems. This eliminates the undesirable effect observed, for example, in \cite{Hauck2025}, where, for a fixed number of element layers in the basis localization, the error increases again after reaching a certain level as the mesh is refined.
    \item In the lowest-order case, the proposed method requires fewer basis functions than \cite{Hauck2025} without notably affecting its convergence or localization behavior. Specifically, while \cite{Hauck2025} uses face integrals as its QOIs, the lowest-order version of the proposed method uses only normal face integrals, reducing the number of basis functions by a factor of $n$ in $n$ dimensions.
    \end{itemize}

The paper is organized as follows. In \cref{sec:modelproblem}, we introduce the heterogeneous Stokes model problem studied in this work. The prototypical multiscale method is presented in \cref{sec:protmethod}. To obtain a practical variant, we show in \cref{sec:expdec} that the method’s basis functions decay exponentially, motivating their localized computation on subdomains. A localized version of the method is then introduced in \cref{sec:locmethod}. The resulting local but still infinite-dimensional subdomain problems are discretized in \cref{sec:finescaledisc} using a fine-scale finite element method. Finally, numerical experiments in \cref{sec:numexp} validate the theoretical findings of the article.

\subsection*{Notation} 
Throughout this work, we use the notation \( a \lesssim b \) (respectively \( b \gtrsim a \)) to indicate that \( a \leq C\, b \) (respectively \( a \geq C\, b \)), where \( C > 0 \) denotes a generic constant independent of the coarse mesh size \( H \), the fine mesh size $h$, the localization parameter \( \ell \), and the oscillations of the PDE solution~\( u \). The constant \(C\) may depend on the mesh regularity, the spatial dimension \(n\), the coefficient bounds \(\nu_{\mathrm{min}}, \nu_{\mathrm{max}}\) and \(\sigma_{\mathrm{min}}, \sigma_{\mathrm{max}}\), and the method order \(m\). The dependence of the constants on $m$ is not tracked explicitly, as no asymptotic behavior with respect to $m$ is considered in this work. Furthermore, we employ standard notation for Sobolev spaces and their norms, denoting by $\|\cdot\|_{k,\Omega}$ the $H^k(\Omega)$-norm, by $|\cdot|_{k,\Omega}$ the highest-order seminorm, and, in the case $k = 0$, we simply write $\|\cdot\|_\Omega$ for the $L^2(\Omega)$-norm.

\section{Model problem}
\label{sec:modelproblem}
This section introduces the weak formulation of the heterogeneous Stokes problem, along with classical results guaranteeing its well-posedness. The formulation is based on the Sobolev space $V \coloneqq (H^1_0(\Omega))^n$, endowed with homogeneous Dirichlet boundary conditions on $\partial \Omega$, and the space $Q \coloneqq \{q \in L^2(\Omega) \with \int_\Omega q \dx = 0\}$, consisting of functions with zero integral mean. In the following, we will always assume the existence of constants $\nu_\mathrm{min}, \nu_\mathrm{max}$ and $\sigma_\mathrm{min}, \sigma_\mathrm{max}$ such that
\begin{equation}
	\label{eq:unifboundcoeff}
	0<\nu_{\min}\le\nu\le\nu_{\max}<\infty,\qquad 0\leq \sigma_{\min}\le\sigma\le\sigma_{\max}<\infty,
\end{equation}
holds almost everywhere in $\Omega$. Denoting by  $(\cdot,\cdot)_\Omega$ the $L^2(\Omega)$-inner product, the problem's bilinear forms $a\colon V\times V \to \mathbb R$ and $b\colon V\times Q \to \mathbb R$ are   defined as
\begin{equation}
	\label{eq:defab}
	a(u,v) \coloneqq (\nu \nabla u ,\nabla v)_{\Omega}  + (\sigma u, v)_{\Omega},\qquad b(u,q) \coloneqq -(q,\div u)_\Omega.
\end{equation}

Given a source term $f \in L^2(\Omega)$, the weak formulation of the considered heterogeneous Stokes problem seeks a pair $(u,p) \in V \times Q$ such that 
\begin{subequations}
	\label{eq:weakstokes}
	\begin{align}
		\qquad\qquad \qquad \qquad &a(u,v)& +\quad &b(v,p) &=\quad &(f,v)_\Omega,&\qquad \qquad \qquad\qquad\label{eq:weakstokes1}\\
		\qquad \qquad \qquad&b(u,q)&   &       &=\quad &0,&\qquad \qquad \qquad\label{eq:weakstokes2}
	\end{align}
\end{subequations}
for all $(v,q) \in V \times Q$.

Using the uniform coefficient bounds~\cref{eq:unifboundcoeff} one can show that the bilinear form~$a$ is coercive and bounded, i.e., there exist constants $c_a, C_a>0$ such that
\begin{equation}
	\label{eq:a}
	|a(v,v)|\geq c_a \|\nabla v\|_\Omega^2,\qquad|a(u,v)| \leq C_a \|\nabla u\|_\Omega\|\nabla v\|_\Omega
\end{equation}
for all functions $u,v \in V$. By the Poincaré--Friedrichs inequality, the seminorm $\|\nabla \cdot\|_\Omega$ is equivalent to the full $(H^1(\Omega))^n$-norm. The constants in \cref{eq:a} can be specified as $c_a = \nu_\mathrm{min}$ and $C_a = \nu_\mathrm{max} + C_\mathrm{PF}^2\sigma_\mathrm{max}$, where $C_\mathrm{PF}>0$ denotes the Poincaré--Friedrichs constant of the domain $\Omega$.  

To establish the well-posedness of problem \cref{eq:weakstokes}, we  need a compatibility condition between the spaces $V$ and $Q$, expressed as the inf--sup condition
\begin{equation}
	\label{eq:infsup}
	\adjustlimits \inf_{q \in Q}\sup_{v \in V} \frac{|b(v,q)|}{\|\nabla v\|_\Omega\|q\|_\Omega}\geq  c_b,
\end{equation}
where $c_b>0$ is typically called the inf--sup constant. This condition is classical and it is typically proved using the so-called Ladyzhenskaya lemma, cf.~\cite{ladyzhenskaia1963}. It states that for any $q \in Q$ there exists $v \in V$ such that
\begin{equation}
	\label{eq:ladlemmaOmega}
	\div v = q,\qquad \|\nabla v\|_\Omega \leq C_\mathrm{L} \|q\|_\Omega,
\end{equation}
which directly implies the inf--sup stability with inf--sup constant $c_b = C_\mathrm{L}^{-1}$. After establishing conditions \cref{eq:a,eq:infsup}, the well-posedness of weak formulation~\cref{eq:weakstokes} can be concluded using classical inf--sup theory; see, e.g.,~\cite{BoffiBrezziFortin2013}.

\section{Prototypical multiscale method}
\label{sec:protmethod}
This section presents a prototypical multiscale method that achieves high-order approximation rates without any pre-asymptotic effects, under minimal structural assumptions on the coefficients. To this end, we introduce a hierarchy of simplicial meshes ${\TH}$ that are geometrically conforming, quasi-uniform, and shape-regular. Each mesh is a finite decomposition of the closure of $\Omega$ into closed elements $T$, which are $n$-dimensional simplices\footnote{The assumption of a simplicial mesh is made only for the simplicity of presentation. Quadrilateral/hexahedral meshes can be used equally well.}. The mesh size, denoted by $H$, is defined as the maximum diameter of the elements in ${\TH}$, i.e., $H \coloneqq \max_{T \in {\TH}} \diam(T)$. For a given polynomial degree $m \in \mathbb{N}_0$, where $\mathbb{N}_0$ denotes the set of natural numbers including zero, we introduce the space $\mathbb{P}^m(T)$ consisting of polynomials of total degree at most $m$ defined on the element $T \in {\TH}$. Piecing together these local spaces in a discontinuous manner gives the space $\mathbb{P}^m({\TH})$ of ${\TH}$-piecewise polynomials of total degree at most $m$. The corresponding $L^2$-orthogonal projection is denoted by $\Pi_H^m\colon L^2(\Omega)\to \mathbb{P}^m({\TH})$. Furthermore, we denote the set of all faces of the mesh ${\TH}$ by $\mathcal{F}_H$ and the subset of interior faces by $\mathcal{F}_H^i$.

The construction of the prototypical high-order LOD method is based on an equivalent reformulation of problem~\cref{eq:weakstokes} using the spaces
\begin{equation}
	\label{eq:defZ}
	Z \coloneqq \big\{ v \in V \with \div v \in \mathbb P^0(\TH)\big\},\qquad Q_H \coloneqq Q \cap \mathbb P^0(\TH),
\end{equation}
where the space $Z$ partially integrates the divergence-free constraint into the velocity space. Thus, the smaller space $Q_H$ is sufficient to enforce that the velocity is divergence-free. The reformulation seeks  $(u,p_H) \in
Z \times Q_H$ such that
\begin{subequations}
	\label{eq:reformulation}
	\begin{align}
		\quad\qquad \qquad \qquad &a(u,v)& +\quad &b(v,p_H) &=\quad &(f,v)_\Omega,&\qquad \qquad \qquad\qquad\label{eq:reformulation1}\\
		\quad \qquad \qquad&b(u,q_H)&   &       &=\quad &0.&\qquad \qquad \qquad
	\end{align}
\end{subequations}
for all $(v,q_H) \in Z\times Q_H$. 
To prove the well-posedness of this reformulated problem, we verify the corresponding inf--sup condition for the bilinear form $b$, which holds with the constant $c_b$ from~\cref{eq:infsup}, thanks again to the Ladyzhenskaya lemma, cf.~\cref{eq:ladlemmaOmega}.
 It is readily seen that the first component of the solution to the reformulated problem coincides with the velocity $u$ from \cref{eq:weakstokes}, while the second component satisfies $p_H = \Pi_H^0 p$, where $p$ is the pressure from \cref{eq:weakstokes}.

\subsection{Quantities of interest}
Following the presentation of the LOD in~\cite{AltHP21}, we introduce quantities of interest (QOIs) that will be preserved by the prototypical method.
 We begin by introducing QOIs associated with faces \( F \in \mathcal{F}_H^i \). Let \( \mathbb{P}^m(F) \) denote the space of polynomials on \( F \) of total degree at most \( m \) and \(\{p_{F,j} \with j = 1, \dots, J\}\) be a basis of \(\mathbb{P}^m(F)\), where \(J \coloneqq \operatorname{dim}(\mathbb{P}^m(F))\). We assume that \(p_{F,1}\equiv 1\)  on~$F$ and that $\int_F p_{F,j} \ds=0$ for indices $j>1$.
  For each face \( F \in \mathcal{F}_H^i \) and index \( j \in \{1, \dots, J\} \), we define the corresponding QOI~as
\begin{align}
	\label{eq:qoiface}
	q_{F,j} \colon Z \to \mathbb R,\qquad  v\mapsto H \int_F (v \cdot n) p_{F,j} \ds,
\end{align}
where $n$ denotes the unit normal vector to $F$, whose direction is fixed once for all.

To obtain a high-order method, it is not sufficient to consider QOIs only on the faces. In addition, also QOIs on the elements must be introduced, similar to the procedure for the MsFEM in {\cite{Feng22}}. A natural choice for these would be the moments against all polynomials in $\mathbb{P}^m (T)$. However, a closer look reveals that one needs to eliminate the moments against the vector-valued polynomials, which are gradients of scalar polynomials of degree at most $m + 1$. The space of such polynomials is in the following denoted as
\begin{align*}
	 \mathbb{G}^m (T) := \{\nabla p \with p \in \mathbb{P}^{m + 1} (T)\}
	\subset (\mathbb{P}^m (T))^n . 
\end{align*}
They must be excluded from the element QOIs, as they conflict with the face QOIs defined in~\cref{eq:qoiface}, as will be explained in more detail in \cref{RemBubblesQm} and the proof of \cref{thm:convergenceprot}. A similar phenomenon also occurs in the Virtual Element Method for the Stokes problem; see~\cite{daVeiga2017}.
At the lowest order \(m = 0\), we have the identity  \((\mathbb{P}^0(T))^n = \mathbb{G}^0(T)\), i.e., no element QOIs need to be considered. For orders \(m \ge 1\), we have the strict subspace relation \(\mathbb{G}^m(T) \subsetneq (\mathbb{P}^m(T))^n\), so the complement of \(\mathbb{G}^m(T)\) in \((\mathbb{P}^m(T))^n\) is nontrivial, and element QOIs need to be defined by choosing a suitable complement of $\mathbb{G}^m(T)$ in $(\mathbb{P}^m(T))^n$. The complement space should satisfy certain conditions, summarized below.

\begin{assumption}[Choice of complement space]
  \label{AssQm}
  Let \(m \ge 1\). Then, for any \({T \in \TH}\), there exists a subspace \(\mathbb{Q}^m(T) \subset (\mathbb{P}^m(T))^n\) such that
  \begin{equation}\label{eq:decomposition}
  	 (\mathbb{P}^m (T))^n =\mathbb{G}^m (T) \oplus \mathbb{Q}^m (T).
  \end{equation}
  Furthermore, for any $p \in (\mathbb{P}^m (T))^n$ the decomposition $p = g + q$, where $g \in \mathbb{G}^m (T)$ and $ q \in \mathbb{Q}^m (T)$ satisfies the following conditions:
  \begin{enumerate}
    \item  $|p|_{\mathbb{P}^m}^2=|g|_{\mathbb{P}^m}^2+|q|_{\mathbb{P}^m}^2$, where $|\cdot|_{\mathbb{P}^m}$ is a seminorm on $\mathbb{P}^m(T)$ defined (independently of $T$) for any $f\in\mathbb{P}^m(T)$ as $|f|_{\mathbb{P}^m}^2\coloneqq \sum_{\alpha:|\alpha|=m}|D^\alpha f|^2$ in the usual notations with multi-indices $\alpha\in\mathbb{N}_0^n$;
    \item the homogeneous parts\footnote{The homogeneous part of degree $m$ of a polynomial is the sum of all terms in the polynomial whose total degree is exactly $m$. }  of degree $m$ of polynomials $g$ and $q$ are determined from the homogeneous part of degree $m$ of polynomial $p$ in a manner independent of the mesh element $T$;
    \item $\|q\|_{T}\lesssim\|p\|_{T}$.
  \end{enumerate}
\end{assumption}

The above conditions on the complement space \(\mathbb{Q}^m(T)\) are assumed to hold throughout the manuscript. A possible construction of such a decomposition in two and three spatial dimensions is given as follows. 
\begin{example}[Construction of the complement space]\label{ExQm2D}
We first consider the two-dimensional case, i.e., \(n = 2\). Let \((x,y)\) denote the coordinate vector and \((x_T, y_T)\) the barycenter of \(T\). To construct the subspaces \(\mathbb G^m(T)\) and \(\mathbb Q^m(T)\), we specify their basis functions. These can be obtained by iterating over pairs $(r,s) \in \mathbb{N}_0^2$ with \(1 \le r+s \le m+1\), and performing the following operations:
	\begin{itemize}
 \item add to the basis of $\mathbb G^m(T)$ the polynomial
 \begin{equation*}
 	\nabla (x - x_T)^r  (y - y_T)^s,
 \end{equation*}
 \item if $r,s>0$ add to the basis of $\mathbb Q^m(T)$ the polynomial
 \begin{equation}
 	\label{eq:basisQ}
 	\left(- r(x - x_T)^{r - 1} (y - y_T)^s,\, 	s(x - x_T)^r (y - y_T)^{s - 1}\right) .
 \end{equation}
	\end{itemize}
Any vector-valued polynomial of the form 
\((\alpha (x - x_T)^{r-1} (y - y_T)^s,\, \beta (x - x_T)^r (y - y_T)^{s-1})\), 
with \(\alpha, \beta \in \mathbb{R}\) and \(r, s>0\) as above, can be represented uniquely as a linear combination of the basis functions in \(\mathbb G^m(T)\) and \(\mathbb Q^m(T)\). Since the span of all such polynomials, plus the polynomials  \((r (x - x_T)^{r-1} ,\, 0)\) and \((0,\, s(y - y_T)^{s-1})\) which are in $\mathbb G^m(T)$, equals \((\mathbb{P}^m(T))^2\), we conclude that our choice of \(\mathbb G^m(T)\) and \(\mathbb Q^m(T)\) satisfies the direct decomposition property~\cref{eq:decomposition}.
Properties (1)--(2) of \cref{AssQm} follow directly from the construction above. Property (3) can be inferred from the following strengthened Cauchy--Schwarz inequality 
$$
(g, q)_{L^2 (T)} \leq \gamma \|g\|_T \|q\|_T,
$$
for all $g\in\mathbb G^m(T)$ and $\ q\in\mathbb Q^m(T)$, with a constant $\gamma<1$ depending only on the regularity of the mesh. This, in turn, can be proved by maximizing the best possible constant \(\gamma\) over all polynomials \(g\) and \(q\) in the corresponding spaces and all triangles \(T\) satisfying the mesh regularity assumption. In this maximization, one can assume, without loss of generality, that \(T\) has its barycenter at the origin and diameter 1, and the $L^2$-norms of $g,q$ are equal to 1. The maximum is thus attained on some element \(\hat{T}\) and some $\hat{g},\hat{q}$, and \(\gamma < 1\), since $\hat{g}$ and $\hat{q}$ are not collinear.

The construction in the three-dimensional case is very similar. The basis functions of \(\mathbb{G}^m(T)\) are of the form 
\(\nabla (x - x_T)^r (y - y_T)^s (z - z_T)^t\) 
for triples \((r,s,t) \in \mathbb{N}_0^3\) with \(1 \le r+s+t \le m+1\). When \(r,s,t > 0\), two basis functions analogous to those in \cref{eq:basisQ} are added to the basis of \(\mathbb{Q}^m(T)\) instead of just one.
\end{example}

Denoting a basis of \(\mathbb{Q}^m(T)\) by \(\{p_{T,k} \with k = 1, \dots, K\}\), where \(K := \dim(\mathbb{Q}^m(T))\), we introduce, for all \(T \in \TH\) and \(k \in \{ 1, \dots, K\}\), the corresponding QOIs as
\begin{align}  
\label{eq:qoielement}
  q_{T, k} : Z \to \mathbb{R}, \qquad v \mapsto \int_T v \cdot p_{T,k} \dx . 
\end{align}

Following the approach in \cite{Hauck2025b}, we handle the QOIs using Lagrange multipliers that belong to the space
\begin{align}
	\label{eq:MH}
	\begin{split}
		M_H \coloneqq \Big\{\bbmu=(\mu,\bmu) \with\ & \mu\colon \Sigma \to\mathbb{R},\ \mu|_F\in\mathbb{P}^m(F)\ \forall F\in\mathcal{F}_H^i, \\
		& \bmu\colon \Omega \to\mathbb{R}^n,\ \bmu|_T\in\mathbb{Q}^m(T)\ \forall T\in\mathcal{T}_H \Big\},
	\end{split}
\end{align}
with the norm
\begin{equation}
	\label{bbmunorm}
	\|\bbmu\|_{M_H}^2 
	\coloneqq 
	H \|\mu\|_\Sigma^2
	+  \|\bmu\|_\Omega^2,
\end{equation}
where \(\Sigma \coloneqq \bigcup_{F \in \mathcal{F}_H^i} F\) denotes the union of all interior faces. We further introduce the bilinear form \(c \colon Z \times M_H \to \mathbb{R}\) as
\begin{equation}
	\label{eq:defc}
	c(v,\bbmu)\coloneqq H \int_{\Sigma} (v \cdot n)  \mu \ds
	+ \int_\Omega v \cdot \bmu \dx,
\end{equation}
with the piecewise defined unit normal
\begin{equation}\label{eq:defn}
	n \colon \Sigma \to \mathbb{R}^n, \qquad n|_F \coloneqq n \quad \forall F \in \mathcal{F}_H^i.
\end{equation}
The bilinear form $c$ encodes the QOIs in the sense that $c(v,\bbmu)=c(w,\bbmu)$ for all $\bbmu\in M_H$ if and only if the functions $v$ and $w$ are indistinguishable with respect to the ensemble of QOIs.  
The scaling with $H$ is included in both the norm on $M_H$ and in the form $c$ to balance the contributions from the faces and the elements.

\subsection{Bubble functions}

Next, we introduce bubble functions, a theoretical tool that will be used repeatedly throughout this manuscript. They play a central role in establishing the inf--sup stability of the bilinear form \(c\), which is crucial for the well-posedness of the proposed method. We consider two types of bubble functions. The first type is used to handle the divergence constraint on each mesh element. The following lemma establishes their existence and summarizes their main properties.

\begin{lemma}[Local Ladyzhenskaya-type bubble]\label{lem:locladymod}
	For any $T \in \mathcal{T}_H$ and any $q \in L^2(T)$ with $\int_T q \, \dx = 0$, there exists $v_q \in (H^1_0(T))^n$ such that $\div v_q = q$, $\int_T v_q \cdot p\dx = 0$ for all $p\in\mathbb{Q}^m(T)$, and the following stability estimate holds:
	\begin{equation}\label{eq:bubbleest3}
		\|\nabla v_q\|_T \lesssim  \|q\|_T.
	\end{equation}
\end{lemma}
\begin{proof}
The proof of this lemma is deferred to \cref{AppendixBubble}.
\end{proof}

The second type of bubble functions corresponds to the constraints associated with the bilinear form $c$ defined in \cref{eq:defc} and the space $M_H$. Their existence and main properties are summarized in the lemma below.

\begin{lemma}[Bubble functions]\label{LemBubbles}
	For any \(F \in \mathcal{F}_H^i\) and any \(g_F \in \mathbb{P}^m(F)\), there exists a face bubble function \(b_F \in Z \cap (H^1_0(\omega_F))^n\) with $\omega_F$ denoting the union of the two mesh elements sharing $F$ such that, for all \(\bbmu \in M_H\),
	\begin{equation}\label{eq:kdfaces}
		c(b_F, \bbmu) = H \int_F g_F \mu \ds,
	\end{equation}
	and
	\begin{equation}\label{eq:bubbleest1}
		\|b_F\|_{\omega_F} \lesssim  H^{1/2} \|g_F\|_F, 
		\qquad 
		\|\nabla b_F\|_{\omega_F} \lesssim  H^{-1/2} \|g_F\|_F.
	\end{equation}
	
	Similarly, for any \(T \in \mathcal{T}_H\) and any \(\mathbf{g}_T \in \mathbb{Q}^m(T)\), there exists an element bubble function \(b_T \in Z \cap (H^1_0(T))^n\) such that, for all \(\bbmu \in M_H\),
	\begin{equation}\label{eq:kdcells}
		c(b_T, \bbmu) = \int_T \mathbf{g}_T \cdot \bmu \dx,
	\end{equation}
	and
	\begin{equation}\label{eq:bubbleest2}
		\|b_T\|_T \lesssim  \|\mathbf{g}_T\|_T, 
		\qquad 
		\|\nabla b_T\|_T \lesssim H^{-1} \|\mathbf{g}_T\|_T.
	\end{equation}
\end{lemma}

\begin{proof}
The proof of this lemma is also deferred to \cref{AppendixBubble}.
\end{proof}

The following remark shows that the construction of the element bubble functions~$b_T$ in the above lemma would not be possible if gradients of scalar polynomials were not excluded from the definition of the element QOIs in \cref{eq:qoielement}.

\begin{remark}[Elimination of gradients in QOIs]\label{RemBubblesQm}
If the definition of the element QOIs in \cref{eq:qoielement} were expanded to include all polynomial moments of degree at most~$m$, that is, if the space $M_H$ were defined by replacing $\mathbb{Q}^m(T)$ in \cref{eq:MH} with $(\mathbb{P}^m(T))^n$, then the construction of the local element bubble function  $b_T$ in \cref{eq:kdcells} would have to accommodate any $\mathbf{g}_T \in (\mathbb{P}^m(T))^n$. In particular, one would need to construct an element bubble function $b_T \in Z \cap (H^1_0(T))^n$ associated with $\mathbf{g}_T = \nabla p$ for some $p \in \mathbb{P}^{m+1}(T)$. Choosing $\boldsymbol{\mu} = \nabla p$ on $T$ in \cref{eq:kdcells} would then yield
\begin{equation}\label{eq:wronginequality}
	\|\nabla p\|_T^2 = \|\mathbf{g}_T\|_T^2 =
	\int_T b_{T} \cdot \nabla p \, dx = - \int_T (\div b_{T}) \, p \, dx = 0,
\end{equation}
since \(\div b_T = 0\) by the divergence theorem, noting that $b_T \in Z \cap (H^1_0(T))^n$. Identity~\cref{eq:wronginequality} is evidently false unless the polynomial $p$ is constant, which reveals an inherent contradiction in such a construction.
\end{remark}

\subsection{Prototypical approximation space}
For defining the prototypical LOD multiscale method, we follow the standard approach in \cite{AltHP21}, decomposing the solution space \(Z\) into a direct sum of two subspaces. The first subspace, often called the fine-scale space, is defined as the intersection of the kernels of the QOIs in \cref{eq:qoiface,eq:qoielement}, which can be expressed using the bilinear form \(c\) introduced in \cref{eq:defc} as
\begin{equation}
	\label{eq:defW}
  W \coloneqq \left\{ v \in Z \with c(v,\bbmu) = 0\ \forall \bbmu \in M_H\right\}.
\end{equation}
This space contains functions that oscillate on scales smaller than \(H\) and cannot be distinguished by the QOIs in~\cref{eq:qoiface,eq:qoielement}.
The second subspace of the decomposition is finite-dimensional and will serve as the approximation space of the prototypical high-order LOD method. It is defined as the orthogonal complement of $W$ with respect to the energy inner product $a$, i.e.,
\begin{equation}
	\label{eq:Zms}
	\tilde Z_H \coloneqq \big\{ u \in Z \with a (u, v) = 0\ \forall  v \in W \big\}.
\end{equation}
Note that, since \(\tilde Z_H\) is constructed as the orthogonal complement of \(W\) with respect to the problem-dependent inner product \(a\), it encodes problem-specific information that allows reliable approximations even on coarse scales. The tildes in the notation of functions and spaces indicate that they are adapted to the problem at hand. The following lemma provides a basis of the space \(\tilde Z_H\).

\begin{lemma}[Prototypical basis]\label{le:protbasis} 
The space \(\tilde{Z}_H\) is of finite dimension \(N := J \cdot \#\mathcal{F}_H^i + K \cdot \#\mathcal{T}_H\), where \(\#(\cdot)\) denotes the number of elements in a set. The basis functions associated with faces, denoted by \(\tilde{\varphi}_{F,j}\) for \(F \in \mathcal{F}_H^i\) and \(j \in \{1, \dots, J\}\), are defined as the unique solution to the problem which seeks \((\tilde{\varphi}_{F,j}, \xi_{F,j}, \bblambda) \in V \times X_H \times M_H\), with \(X_H \coloneqq \{q \in Q \with \Pi_H^0 q = 0\}\), such that
\begin{subequations}\label{pbphiE}
	 \begin{align}
		&  \quad \qquad a (\tilde{\varphi}_{F, j}, v) & + \quad & b (v,
		\xi_{F, j}) & + \quad & c (v, \bblambda) & = \quad & 0, \quad  \quad &
		\label{eq:LODbasis1}\\
		& \quad \qquad b (\tilde{\varphi}_{F, j}, \chi) &  &  &  &  & =
		\quad & 0, \quad \qquad &  \label{eq:LODbasis2}\\
		& \quad \qquad c (\tilde{\varphi}_{F, j}, \bbmu) &  &  &  &  & =
        \quad & H \tint_F p_{F,j}\mu \ds \quad &  \label{eq:LODbasis3}
	\end{align}
\end{subequations}  
  for all $(v, \chi, \bbmu) \in V \times X_H \times M_H$.
  
  Moreover, the basis functions associated with elements, denoted by \(\tilde{\varphi}_{T,k}\) for all \(T \in \mathcal{T}_H\) and \(k \in \{1, \dots, K\}\), are defined as the solution to a problem analogously to~\cref{pbphiE}, with the indices \(F\) and \(j\) replaced by \(T\) and \(k\), respectively, and the right-hand side of the last equation replaced by $\int_T p_{T,k} \cdot \bmu \dx$.
\end{lemma}
\begin{proof}
We interpret (\ref{pbphiE}) as a standard saddle-point problem, combining the Lagrange multiplier spaces $X_H$ and $M_H$ into the product space $X_H\times M_H$ equipped with the norm  $\|(q,\bbmu)\|^2 \coloneqq \|q\|_\Omega^2+H^2\|\bbmu\|_{M_H}^2$ for $(q,\bbmu) \in X_H\times M$. The well-posedness of this problem then follows from the inf--sup condition:
\begin{equation}\label{infsupcbb}
		\adjustlimits  \inf_{(q,\bbmu) \in X_H\times M_H} \sup_{v \in V} \, \frac{b(v,q)+c(v,\bbmu)}{\|\nabla v\|_\Omega \|(q,\bbmu)\|}\gtrsim 1.
\end{equation}
To prove this, take any $(q,\bbmu) \in X_H \times M_H$, and, thanks to \cref{lem:locladymod}, introduce $v_{q,T} \in (H^1_0(T))^n$ for each $T \in \TH$ such that $b(v_{q,T},q) = \|q\|_T^2$, $c(v_{q,T},\bbmu) = 0$ for all $\bbmu \in M_H$, and $\|\nabla v_{q,T}\|_T \lesssim \|q\|_T$. Similarly, for each $F \in \mathcal{F}_H^i$ (resp. each $T \in \TH$), we introduce, thanks to \cref{LemBubbles}, $v_{\mu,F} \in Z \cap (H^1_0(\omega_F))^n$ (resp. $v_{\bmu,T} \in Z \cap (H^1_0(T))^n$) such that $c(v_{\mu,F},\bbmu) = H \|\mu\|_F^2$ with $\|\nabla v_{\mu,F}\|_{\omega_F} \lesssim H^{-1/2} \|\mu\|_F$ (resp. $c(v_{\bmu,T},\bbmu) = \|\bmu\|_T^2$ with $\|\nabla v_{\bmu,T}\|_T \lesssim H^{-1} \|\bmu\|_T$). Setting
\begin{equation*}
	v\coloneqq \sum_{T\in\TH}v_{q,T}+H^2\sum_{F\in\mathcal{F}_H^i}v_{\mu,F}+H^2\sum_{T\in\TH}v_{\bmu,T},
\end{equation*}
we obtain $b(v,q)+c(v,\bbmu)=\|(q,\bbmu)\|^2$ as well as
\begin{align*}
	\|\nabla v\|_\Omega^2
	&\lesssim \sum_{T\in\TH}\|\nabla v_{q,T}\|_T^2+H^4\sum_{F\in\mathcal{F}_H^i}\|\nabla v_{\mu,F}\|_{\omega_F}^2+H^4\sum_{T\in\TH}\|\nabla v_{\bmu,T}\|_T^2	\lesssim \|(q,\bbmu)\|^2,
\end{align*}
which proves \cref{infsupcbb}, the desired inf--sup condition.

By construction, the functions in the set $\{\tilde \varphi_{F,j}\}_{F,j} \cup \{\tilde \varphi_{T,k}\}_{T,k}$ belong to $\tilde Z_H$. 
To verify that they in fact form a basis of this space, we consider any $u\in \tilde Z_H$ and observe that there exists a unique linear combination of $\{\tilde \varphi_{F,j}\}\cup\{\tilde \varphi_{T,k}\}$, say $w\in \tilde Z_H$, such that $c(w,\bbmu)=c(u,\bbmu)$ for all $\bbmu\in M_H$. The well-posedness of the saddle-point problem for $u-w$, as above, then implies that $u=w$. 
\end{proof}

The functions in the fine-scale space $W$ satisfy an element-local Poincaré-like property, as stated in the following lemma.
\begin{lemma}[Local Poincare-type inequality]\label{lem:poincare}
	For all $T \in \TH$ and $v \in (H^1(T))^n$ with $\tint_F v\cdot n \ds  = 0$ for all faces $F \subset \partial T$, it holds that
	\begin{equation}
		\label{eq:Poincarelocal}
		\|v\|_T \lesssim H \|\nabla v\|_T.
	\end{equation}
\end{lemma}
\begin{proof}
    The proof of this Lemma is deferred to \cref{AppendixBubble}. 
\end{proof}
We emphasize that, unlike in \cite{Hauck2025}, where face integrals of vector-valued functions $v$ were used as QOIs, we now consider only (suitably weighted) normal integrals over faces in the definition of the face QOIs in~\cref{eq:qoiface}.  

\subsection{Prototypical method}
Having introduced the prototypical approximation space $\tilde Z_H$, we define the prototypical method by replacing $Z$ in \cref{eq:reformulation} with $\tilde Z_H$. Specifically, it seeks $(\tilde u_H, \tilde p_H) \in \tilde Z_H \times Q_H$ such that
\begin{subequations}
	\label{LODid2}
	\begin{align}
		\qquad \qquad \qquad &a(\tilde u_H,\tilde v_H)& +\quad &b(\tilde v_H,\tilde p_H) &=\quad &(f,\tilde v_H)_\Omega,&\qquad \qquad \qquad\label{eq:LODdivfree1}\\
		\qquad \qquad \qquad&b(\tilde u_H,\tilde q_H)&   &       &=\quad &0&\qquad \qquad \qquad\label{eq:LODdivfree2}
	\end{align}
\end{subequations}
for all $(\tilde v_H,\tilde q_H) \in \tilde Z_H\times Q_H$. 

To characterize the solution of the prototypical method and as a tool in its analysis,  we introduce the $a$-orthogonal projection operator $\mathcal{R} : Z \to \tilde Z_H$. Given any $v \in Z$, its projection $\mathcal{R}v$ is defined as the unique element of $\tilde Z_H$ satisfying
\[
a(v - \mathcal{R} v, w) = 0\quad \forall w \in \tilde Z_H.
\]
The $a$-orthogonality immediately yields the continuity of $\mathcal{R}$. Indeed, for all $v \in Z$,
\begin{equation} \label{eq:Rcont} \|\nabla \mathcal{R} v\|_\Omega \le \sqrt{C_a / c_a} \, \|\nabla v\|_\Omega, \end{equation}
where $c_a$ and $C_a$ denote the constants from~\cref{eq:a}.
Since $\tilde Z_H$ is the $a$-orthogonal complement of the fine-scale subspace $W$, we also have, for any $v \in Z$,
\[
c(v-\mathcal{R} v, \bbmu) = 0\quad \forall \bbmu \in M_H
\]
as $v-\mathcal{R} v\in W$. Taking the face-based component of $\bbmu$ above as piecewise constant shows that the operator $\mathcal{R}$ preserves face normal fluxes, that is $\int_F(v-\mathcal{R} v)\cdot n\ds=0$ on all faces $F \in \mathcal F_H^i$. Thus, using the divergence theorem and noting that functions in $Z$ have piecewise constant divergence (cf.~\cref{eq:defZ}), we obtain, for all $v \in Z$,
\begin{equation}
	\label{eq:identityIH}
	(\div \mathcal{R} v)|_T = (\div v)|_T\quad \forall T \in \mathcal{T}_H\,.
\end{equation}

The following theorem provides a convergence result for the prototypical method, valid under minimal structural assumptions on the coefficients.

\begin{theorem}[Prototypical method]\label{thm:convergenceprot}
	The prototypical multiscale method~\cref{LODid2} is well-posed, and its solution is given by $(\tilde u_H,\tilde p_H)=(\mathcal R u,\Pi_H^0p)$, where $(u,p)$ solves problem~\cref{eq:weakstokes}.
Moreover, for any \( f \in H^{m+1}(\Omega) \), we have the error estimates
	\begin{align}
		\|\nabla (u - \tilde u_H)\|_\Omega &\lesssim H^{m+2} |f|_{m+1,\Omega},\label{eq:errestH1}\\
		\|u - \tilde u_H\|_\Omega &\lesssim H^{m+3} |f|_{m+1,\Omega}.\label{eq:errestL2}
	\end{align}
\end{theorem}
\begin{proof}
	
	First, we prove the  inf--sup condition
	\begin{equation}\label{eq:infsupLOD}
		\adjustlimits \inf_{\tilde q_H \in Q_H}\sup_{\tilde v_H \in \tilde Z_H} \frac{|b(\tilde v_H,\tilde q_H)|}{\|\nabla \tilde v_H\|_\Omega\|\tilde q_H\|_\Omega} \gtrsim 1,
	\end{equation}
	which implies the well-posedness of problem~\cref{LODid2}.
Given any $\tilde q_H \in Q_H$, let $v \in V$ satisfy $\div v = \tilde q_H$ and 
$\|\nabla v\|_\Omega \lesssim  \|\tilde q_H\|_\Omega$, cf.~\cref{eq:ladlemmaOmega}, 
and define $\tilde v_H \coloneqq \mathcal R v$.  
Using~\cref{eq:identityIH}, we have 
$|b(\tilde v_H, \tilde q_H)| = |b(v, \tilde q_H)| = \|\tilde q_H\|_\Omega^2$, 
and from~\cref{eq:Rcont} and the choice of $v$, 
$\|\nabla \tilde v_H\|_\Omega \lesssim  \|\nabla v\|_\Omega \lesssim  \|\tilde q_H\|_\Omega$.  inf--sup condition \cref{eq:infsupLOD} follows after combing these results.
	Since $(u,\Pi_H^0p)$ solves~\cref{eq:reformulation}, $\mathcal R\colon Z \to \tilde Z_H$ is the $a$-orthogonal projection, and $\div\tilde{u}_H=0$ by~\cref{eq:identityIH}, we observe that $(\mathcal Ru,\Pi_H^0p)$ solves \cref{LODid2}. The uniqueness of the solution to \cref{LODid2} then implies that $\tilde u_H = \mathcal Ru$ and $\tilde p_H = \Pi_H^0p$.

Next, we show the convergence of the prototypical method. Let us denote the error by 
$e \coloneqq u - \tilde u_H$. It holds that $e \in W$, and we have for the prototypical method that 
$\div \tilde u_H = 0$.  
Using the coercivity of $a$ (cf.~\cref{eq:a}), the orthogonality $a(\tilde u_H, e) = 0$, 
\cref{eq:reformulation1} with $e$ as the test function, and the fact that $b(e, p_H) = 0$, we obtain that
	\begin{equation}
		\label{eq:errestprot}
		c_a\|\nabla e\|_\Omega^2
		\leq a (e, e) = a (u, e) = (f, e)_{\Omega}
		= (f-\Pi_H^m f,e)_\Omega + (\Pi_H^m f,e)_\Omega.
	\end{equation}
	The first term on the right-hand side above can be bounded using a classical approximation result for the $L^2$-projection onto piecewise polynomials, cf. \cite[Lem.~1.58]{PiE12}, and the Poincaré-type inequality from \cref{lem:poincare} as
	\begin{equation}\label{eq:est1}
		(f - \Pi_H^m f, e)_{\Omega} \leq \|f - \Pi_H^m f\|_{\Omega} \|e\|_{\Omega}
		\lesssim  H^{m+2} |f|_{m+1,\Omega}  \| \nabla e\|_{\Omega}.
	\end{equation}
	
	To estimate the second term on the right-hand side of \cref{eq:errestprot}, we abbreviate 
	$f_T \coloneqq (\Pi_H^m f)|_T$ for each $T \in \mathcal{T}_H$, and decompose 
	$f_T = g_T + q_T$ with $g_T \in \mathbb{G}^m(T)$ and $q_T \in \mathbb{Q}^m(T)$ as in \cref{AssQm}. 
	By construction, $g_T = \nabla \phi_T$ for some $\phi_T \in \mathbb{P}^{m+1}(T)$.  
	Since $e$ is $L^2(T)$-orthogonal to $\mathbb{Q}^m(T)$ on all elements and $\div e = 0$, the divergence theorem, applied locally on $T$, gives 
	\begin{equation}\label{trickyCR}
		(\Pi_H^m f, e)_{\Omega} = \sum_{T \in \mathcal{T}_H} (g_T + q_T, e)_T 
		= \sum_{T \in \mathcal{T}_H} (\nabla \phi_T, e)_T
		= \sum_{F \in \mathcal{F}_H^i} ([\phi]_F, e \cdot n)_F,
	\end{equation}
	where $[\phi]_F = \phi_T - \phi_{T'}$ denotes the jump across face $F$ shared by elements $T$ and~$T'$, ordered consistently with the normal $n$ on $F$. 

Since \(e \cdot n\) is \(L^2(F)\)-orthogonal to \(\mathbb{P}^m(F)\) for all faces \(F \in \mathcal{F}_H^i\), we can subtract from  \([\phi]_F\) any \(\lambda_F \in \mathbb{P}^m(F)\) in the right-hand side of equation~\cref{trickyCR}. Taking $\lambda_F$ as the \(L^2(F)\)-projection of \([\phi]_F\) to $\mathbb{P}^m(F)$, and using standard approximation properties of the \(L^2\)-projection, cf.~\cite[Lem.~1.58]{PiE12}, we can estimate \cref{trickyCR} as
	\begin{multline}
		\label{eq:est12}
			\sum_{F \in \mathcal{F}_H^i} ([\phi]_F, e \cdot n)_F
			= \sum_{F \in \mathcal{F}_H^i} ([\phi]_F - \lambda_F, e \cdot n)_F \\
			\lesssim H^{m+1} \sum_{F \in \mathcal{F}_H^i} |[\phi]_F|_{m+1,F} \| e \cdot n \|_F 
			\le H^{m+1} \sum_{F \in \mathcal{F}_H^i} {|F|}^{\frac 12}|g_T - g_{T'}|_{\mathbb{P}^m} \| e \|_F\,,
	\end{multline}
where $|\cdot|_{\mathbb{P}^m}$ is the seminorm on $\mathbb{P}^m(T)$ defined in property (1) of \cref{AssQm}. We also have
\begin{equation*}
	|f_T - f_{T'}|_{\mathbb{P}^m}^2 = |g_T - g_{T'}|_{\mathbb{P}^m}^2 + |q_T - q_{T'}|_{\mathbb{P}^m}^2\,,
\end{equation*}
by properties (1)--(2) of the same assumption.
Thus,
\[
{|F|}^{\frac 12}|g_T - g_{T'}|_{\mathbb{P}^m}
\le {|F|}^{\frac 12}|f_T - f_{T'}|_{\mathbb{P}^m}
= |[f - \Pi_H^m f]|_{H^{m,\textrm{full}}(F)}
\lesssim \sqrt{H}\, |f|_{m+1,\omega_F},
\]
where $\omega_F=T\cup T'$ and the seminorm $|\cdot|_{H^{m,\textrm{full}}(F)}$ includes all partial derivatives of order $m$, not only those tangential to $F$. The last inequality follows from a standard approximation result in \(L^2(T)\) and  \(L^2(T')\) and the trace inequality.
Finally, inserting the latter estimate into~\cref{eq:est12}, and using the bound \(\| e \|_F \lesssim \sqrt{H} \| \nabla e \|_{\omega_F}\), which can be derived from \cref{lem:poincare} and a standard trace inequality, we obtain 	\begin{equation}\label{eq:est2}
		(\Pi_H^m f, e)_{\Omega} \lesssim  H^{m+2} \sum_{F \in \mathcal F_H^i} |f|_{m+1,\omega_F} \|\nabla e\|_{\omega_F} 
		\lesssim  H^{m+2} |f|_{m+1,\Omega} \|\nabla e\|_{\Omega}.
	\end{equation}
The desired $H^1$-estimate \cref{eq:errestH1} follows by inserting \cref{eq:est1,eq:est2} into~\cref{eq:errestprot}.

The $L^2$-estimate \cref{eq:errestL2} is obtained by applying \cref{lem:poincare} once again. 
\end{proof}

\section{Exponential decay and localization}
\label{sec:expdec}
We emphasize that the prototypical LOD basis functions defined in \cref{pbphiE} are globally supported. Consequently, computing them would require solving global problems, which we consider infeasible in practice. In this section, we show that the prototypical LOD basis functions decay exponentially, which motivates their approximation by locally computable counterparts. A practical multiscale method based on such local approximations is presented in \cref{sec:locmethod}.

The localization strategy used below was first introduced in~\cite{Hauck2025b} and is highly flexible, making it suitable for applications such as Stokes problems that go beyond standard elliptic diffusion-type problems.  Its key idea is to decompose
\begin{equation}\label{defK}
	\mathcal{R} = \mathcal{I}_H + \mathcal{K},
\end{equation}
where $\mathcal{K} \colon V \to V$ is characterized below, and $\mathcal{I}_H \colon V \to V_H$ with $V_H \subset V$ is a quasi-interpolation operator onto the standard conforming finite element space of (vector-valued) piecewise affine functions on $\mathcal{T}_H$.
The operator $\mathcal{I}_H$ is assumed to depend only on the normal face integrals $\int_F v \cdot n \,\ds$ for each $F \in \mathcal{F}_H^i$ of an input~$v$, and to satisfy the standard approximation and stability properties:
\begin{equation}
	\label{intEst}
	\|v - \mathcal{I}_H v\|_T \lesssim H \| \nabla v \|_{\mathsf N(T)},
	\quad\|\mathcal I_H v\|_T \lesssim \|v\|_{\Nb(T)},\quad 
	\| \nabla \mathcal{I}_H v \|_T \lesssim \| \nabla v \|_{\mathsf N(T)},
\end{equation}
for all functions \(v \in V\) and for any element \(T \in \mathcal{T}_H\). Here, for a union of coarse elements \(S\), \(\mathsf{N}(S)\) denotes the first-order patch consisting of all coarse elements that share at least one node with an element in \(S\). We emphasize  that \(V_H\) is fixed as the first-order finite element space, independent of the polynomial order \(m\).

\begin{remark}[Possible construction of $\mathcal I_H$]
A possible definition of $\mathcal I_H$ in the 
two-dimensional case is given by prescribing its nodal values at interior nodes $z$ as
	\begin{equation*}
		(\mathcal I_H v)(z)
		\coloneqq  
		\begin{bmatrix}
			n_{F_1}^1 & n_{F_1}^2\\
			n_{F_2}^1 & n_{F_2}^2
		\end{bmatrix}^{-1}
		\begin{bmatrix}
			|F_1|^{-1}\int_{F_1} v \cdot n_{F_1}\, \mathrm ds \\
			|F_2|^{-1}\int_{F_2} v \cdot n_{F_2}\, \mathrm ds
		\end{bmatrix},
	\end{equation*}
	where $F_1$ and $F_2$ are any two faces adjacent to $z$ whose normal vectors 
	$n_{F_1}$ and~$n_{F_2}$ (with superscripts denoting their components) are linearly independent.
	For boundary nodes $z$ we set $(\mathcal I_H v)(z) \coloneqq 0$. This definition 
	can be extended to the three-dimensional case by selecting, for all interior nodes $z$, three faces adjacent to $z$ with linearly independent normal vectors. The stability and approximation properties in \cref{intEst} can be verified for this operator following the approach of \cite[Ch.~1.6]{ErG04}.
\end{remark}

The operator \( \mathcal{K} \colon V \to V \), as introduced in \cref{defK}, is characterized for each \( v \in V \) as the unique solution $(\mathcal{K}v,\xi,\bblambda) \in V\times X_H\times M_H$ satisfying
\begin{align*}
	&\qquad \qquad a (\mathcal{K}v, w)& +\quad  &b(w,\xi)          & +\quad  &c(w,\bblambda) & =\quad  &-a (\mathcal{I}_Hv, w),\qquad \quad\quad &\\
	& \qquad \qquad  b(\mathcal{K}v,\chi)                &   &&   &             & =\quad  &-b(\mathcal{I}_Hv,\chi),\qquad\quad\qquad &\\
	&\qquad\qquad c(\mathcal{K}v,\bbmu)  &   &                  &   &             & =\quad  &c(v-\mathcal{I}_Hv,\bbmu)\qquad \quad\quad&
\end{align*}
for all $(w,\chi,\bbmu) \in V\times X_H\times M_H$. The operator \(\mathcal{K}\) can be represented as the following sum of localizable element contributions:
\begin{equation}
	\label{eq:Ksum}
	\mathcal{K} = \sum_{T \in \mathcal{T}_H} \mathcal{K}_T,
\end{equation}
where \(\mathcal{K}_T \colon V \to V\) is defined, for any \(T \in \mathcal{T}_H\) and all \(v \in V\), as the unique solution to the problem, which seeks  \((\mathcal{K}_T v, \xi_T, \bblambda_T) \in V \times X_H \times M_H\) such that
\begin{subequations}\label{eq:KT}
	\begin{align}
		&\qquad  a(\mathcal{K}_Tv, w)& +\quad  &b(w,\xi_T)          & +\quad  &c(w,\bblambda_T) & =\quad  &-a_T(\mathcal{I}_Hv, w),\qquad  &\label{eq:KTa}\\
		& \qquad  b(\mathcal{K}_Tv,\chi)                &   &&   &             & =\quad  &-b_T(\mathcal{I}_Hv,\chi),\qquad &\label{eq:KTb}\\
		&\qquad c(\mathcal{K}_Tv,\bbmu)  &   &                  &   &             & =\quad  &c_T(v-\mathcal{I}_Hv,\bbmu)\qquad&\label{eq:KTc}
	\end{align}
\end{subequations}
for all $(w,\chi,\bbmu) \in V\times X_H\times M_H$. Here, \(a_T\) and \(b_T\) denote the restrictions of \(a\) and~\(b\) to \(T\) (i.e., integrals over \(\Omega\) are replaced by integrals over \(T\)), and \(c_T\) is defined as
\begin{equation}\label{defcT}
    c_T(v,\bbmu) \coloneqq
         H \int_{\partial T} \kappa_{T} (v \cdot n)\mu \ds +
    \int_{T} v \cdot \bmu\dx,
\end{equation}
where the function \(\kappa_T \colon \Sigma \to \mathbb{R}\) is piecewise defined, with \(\kappa_T|_F \ge 0\) nonzero only for faces satisfying \(\omega_F \supset T\), and such that \(\sum_{T \in \omega_F} \kappa_T|_F = 1\) for all \(F \in \mathcal{F}_H^i\). The unit normal \(n \colon \Sigma \to \mathbb{R}^n\) is defined in \cref{eq:defn}.

To quantify the decay of the operators $\mathcal K_T$, we introduce the notion of patches with respect to the mesh \(\mathcal{T}_H\). Given an localization parameter \(\ell \in \mathbb{N}\), the \(\ell\)-th order patch of a set of elements \(S \subset \mathcal{T}_H\) is defined by
\begin{equation}
	\Nb^1(S) \coloneqq \Nb(S), \qquad 
	\Nb^\ell(S) \coloneqq \Nb^1(\Nb^{\ell-1}(S)), \quad \ell \ge 2,
\end{equation}
where we recall that \(\Nb(S)\) denotes the first-order patch of elements around \(S\).
The following theorem shows that the operators \(\mathcal{K}_T\) exhibit an exponential decay.

\begin{theorem}[Exponential decay]
	\label{lem:expdecay}
	There exists a constant $c>0$ independent of $H,\ell,T$, such that, for any \(v \in V\) and $\ell \in \mathbb N$, it holds that
	\[
	\|\nabla \mathcal{K}_T v\|_{\Omega \setminus \Nb^{\ell}(T)} \lesssim \exp(- \mathfrak{c} \ell) 
	\, \|\nabla \mathcal{K}_T v\|_{\Omega}.
	\]
\end{theorem}

\begin{proof}
	 This proof follows the arguments used in the proof of \cite[Thm.~4.1]{Hauck2025} for the lowest-order case. Let \(\eta \in W^{1,\infty}(\Omega)\) be a cut-off function such that:
  \begin{align*}
    \left\{ \begin{aligned}
      \eta & \equiv 0 & \quad & \text{in } \Nb^{\ell - 1} (T),\\
      \eta & \equiv 1 & \quad & \text{in } \Omega \setminus \Nb^{\ell} (T),\\
      0 & \leq \eta \leq 1 & \quad & \text{in } R \coloneqq \overline{ \Nb^{\ell} (T)
      \setminus \Nb^{\ell - 1} (T)},
    \end{aligned} \right.  \label{eq:eta}
  \end{align*}
satisfying the Lipschitz bound \(\sup_{\Omega} |\nabla \eta| \lesssim H^{-1}\), where we used \(\overline{\cdot}\) to denote the closure of a set. The notation \(\mathring{\cdot}\) will be used to denote the interior of a set.
  
 Abbreviating $\psi \coloneqq \mathcal{K}_T v$ and choosing $\eta \psi$ as a test function in~\cref{eq:KTa}, we obtain
 \begin{equation}
 	\label{eq:idabc}
 	a(\psi, \eta \psi)
 	= -\, b(\eta \psi, \xi)
 	- c(\eta \psi, {\bblambda})
 	- a_T(\mathcal{I}_H v, \eta \psi),
 \end{equation}
 where, for simplicity, we have omitted the subscript~$T$ on $\xi_T$ and ${\bblambda}_T$. 
Note that $a_T(\mathcal{I}_H v, \eta \psi) = 0$. 
Moreover, $\operatorname{supp}(\eta \psi) \subset (\Omega \setminus \Nb^{\ell}(T)) \cup R$, 
$\nabla\!\cdot(\eta \psi)$ is piecewise constant on $\Omega \setminus R$, 
$\int_{E}\eta \psi\cdot n\mu\ds = 0$ for all $\mu\in\mathbb{P}^m(F)$ on any face~$E$ not contained in~$\mathring{R}$, 
and $\int_{K}\eta \psi\cdot\bmu\dx = 0$ for all $\mu\in\mathbb{Q}^m(K)$ on any element~$K$ not contained in~$R$. 
Using these properties, together with bound~\cref{eq:a}, we can rewrite~\cref{eq:idabc} as
 \begin{equation}
 	\label{eq:identsubd}
 	c_a \, \|\nabla \psi\|_{\Omega \setminus \Nb^{\ell}(T)}^2
 	\le a_{\Omega \setminus \Nb^{\ell}(T)}(\psi, \psi)
 	= -\,\underbrace{a_{R}(\psi, \eta \psi)}_{\eqqcolon \Xi_1}
 	-\,\underbrace{b_{R}(\eta \psi, \xi)}_{\eqqcolon \Xi_2}
 	-\,\underbrace{c_{R}(\eta \psi, {\bblambda})}_{\eqqcolon \Xi_3}.
 \end{equation}
Here, the restricted forms $a$, $b$, and $c$ are obtained by restricting the integrals in~\cref{eq:defab} to the corresponding subdomain, and by restricting the sums in \cref{eq:defc} to the faces in $\mathring{R}$ and the elements in $R$ for $c_R$, which should not be confused with \cref{defcT}.
   
To estimate \(\Xi_1\), we note that the Poincaré-type inequality from \cref{lem:poincare} can be applied locally to the function \(\psi\) on all elements \(K \subset R\) for \(\ell \ge 2\), since, by~\cref{eq:KTc}, \(\int_E \psi\cdot n \ds = 0\) for all faces \(E \subset \partial K\). Using this property, together with the definition of bilinear form $a$ in \cref{eq:defab} and the \(L^\infty\)- and Lipschitz bounds of \(\eta\), we obtain the following estimate for \(\ell \ge 2\):
  \begin{align*}
  \Xi_1 &
  \lesssim \|\nabla \psi\|_{R} \left( \|
    \nabla \psi\|_R + H^{- 1} \|\psi\|_R \right) + \|\psi\|^2_R\lesssim \|\nabla \psi\|_{
    R}^2.
  \end{align*}
  
For the term \(\Xi_2\), we again apply the \(L^\infty\)- and Lipschitz bounds of \(\eta\) together with the local Poincaré-type inequality from \cref{lem:poincare} for \(\psi\) to obtain that
  \begin{equation}\label{eq:Xi1}
  	\Xi_2  \leq \|\div (\eta \psi) \|_{R}\|\xi\|_{R}
  	\lesssim \|\nabla \psi\|_{R}\|\xi\|_{R}.
  \end{equation}
We continue the previous estimate by deriving a bound for \(\|\xi\|_{K}\) on any element \(K \subset R\). To this end, we test \cref{eq:KTa} with \(v_\xi \in (H^1_0(K))^n\) chosen such that \(\div v_\xi = \xi\) holds locally in \(K\), \(q_{K,l}(v_\xi) = 0\) for all indices \(l\), and \(\|\nabla v_\xi\|_K \lesssim \|\xi\|_K\), where the existence of such a function \(v_\xi\) is guaranteed by \cref{lem:locladymod}. This yields
  \begin{align*}
  	\|\xi\|_K ^2 = -b(v_\xi,\xi) = a(\psi,v_\xi) + c(v_\xi,\bblambda) + a_T(\mathcal I_H v,v_\xi)\lesssim  \|\nabla\psi\|_K\|\xi\|_K,
  \end{align*}
using that \(c(v_\xi, \bblambda) = 0\) by the construction of \(v_\xi\), and that \(a_T(\mathcal I_H v, v_\xi) = 0\) since \(T\) is not contained in \(R\). Summing the above bound over all \(K \subset R\) yields an estimate for \(\|\xi\|_{R}\), which can be inserted into \cref{eq:Xi1} to conclude that \(\Xi_2 \lesssim \|\nabla {\psi}\|_{R}^2.\)

The term $\Xi_3$ can be estimated as
  \begin{align}\label{eq:Xi3}
  	    \Xi_3 &\lesssim H\sum_{E \subset \mathring{R}}\|\eta \psi\|_{E}\|\lambda\|_{E} + \sum_{K \subset R\vphantom{\mathring{R}}}\|\eta \psi \|_K \|\blambda\|_K.
  \end{align}
To estimate \(\|\lambda\|_E\) for any face \(E \subset \mathring{R}\), we consider a bubble function \(b_E \in Z \cap (H^1_0(\omega_E))^n\) satisfying \(c(b_E, \bbmu) = -H\int_E \lambda \mu \ds\) for all \(\bbmu \in M_H\), whose existence is guaranteed by \cref{LemBubbles}. Testing \cref{eq:KTa} with $b_E$  yields
  \begin{align*}
  	H\|\lambda\|_E^2 = -c(b_E,\bblambda) = a(\psi,b_E) + b(b_E,\xi) + a_T(\mathcal I_H v,b_E) \lesssim H^{-1/2}\|\nabla \psi\|_{\omega_E}\|\lambda\|_{E},
  \end{align*}
where we have used that \(b(b_E, p) = 0\), since \(b_E\) has piecewise constant divergence and \(\xi\) has zero element averages, and that \(a_T(\mathcal I_H v, b_E) = 0\) because \(\omega_E\), which is contained in $R$, has only a trivial intersection with \(T\).

To estimate \(\|\blambda\|_K\) for any element \(K \subset R\), we consider a bubble function \(b_K \in Z \cap (H^1_0(K))^n\) that satisfies \(c(b_K, \bbmu) = -\int_K \blambda \cdot \bmu \, \dx\) for all \(\bbmu \in M_H\). Testing~\cref{eq:KTa} with \(b_K\), we obtain, similarly to the previous estimate, that
    \begin{align*}
    	\label{eq:xi3expdec}
  	\|\blambda\|_K^2 = -c(b_K,\bblambda) = a(\psi,b_K) + b(b_K,\xi) + a_T(\mathcal I_H v,b_K) \lesssim H^{-1}\|\nabla \psi\|_{K}\|\blambda\|_{K},
  \end{align*}
noting that \(b_K \in Z \cap (H^1_0(K))^n\) already implies that \(b_K\) is divergence-free.

We can now combine the previous two estimates for \(\|\lambda\|_E\) and \(\|\blambda\|_K\) to continue estimate \cref{eq:Xi3}. Using a standard trace inequality and the Poincaré-type inequality from \cref{lem:poincare} applied to \(\psi\), we then obtain, for \(\ell \ge 2\), that
  \begin{align*}
  	\Xi_3 \lesssim \sum_{E \subset \mathring{R}}\big(H^{-1/2}\|\psi\|_{\omega_E} + H^{1/2}\|\nabla \psi\|_{\omega_E}\big)H^{-1/2}\|\nabla \psi\|_{\omega_E} + \|\nabla \psi\|_R^2 \lesssim \|\nabla \psi\|_R^2.
  \end{align*}

Inserting the above estimates for \(\Xi_1\)--\(\Xi_3\) into \cref{eq:identsubd} gives
\begin{equation*}
	\|\nabla \psi\|_{\Omega \setminus N^{\ell}(T)}^2 \le C \|\nabla \psi\|_{R}^2
	= C \Big(\|\nabla \psi\|_{\Omega \setminus N^{\ell - 1}(T)}^2 - \|\nabla \psi\|_{\Omega \setminus N^{\ell}(T)}^2\Big),
\end{equation*}
where the constant \(C > 0\) is independent of \(H\), \(\ell\), and \(T\). This immediately yields
\[
\|\nabla \psi\|_{\Omega \setminus N^{\ell}(T)} \le \left(\frac{C}{1 + C}\right)^{1/2} \, \|\nabla \psi\|_{\Omega \setminus N^{\ell - 1}(T)},
\]
and, after iterating, it follows that
\[
\|\nabla \psi\|_{\Omega \setminus N^{\ell}(T)} \le \left( \frac{C}{1 + C} \right)^{\ell/2} \|\nabla \psi\|_\Omega
= \exp(- \mathfrak{c} \ell) \, \|\nabla \psi\|_\Omega,
\]
with \( \mathfrak{c} := \frac{1}{2} \log \frac{1 + C}{C} > 0\).
\end{proof}
The exponential decay result from the previous theorem motivates the localization of the operators \(\mathcal{K}_T\) to \(\ell\)-th order patches around \(T\). To this end, we introduce, for all \(T \in \mathcal{T}_H\) and a given $\ell \in \mathbb N$, local versions of the spaces \(V\), \(X_H\), and \(M_H\) as
\begin{align*}
	V_{T}^\ell &\coloneqq \{v \in V \with \operatorname{supp}(v) \subset \Nb^\ell(T)\},\\
	X_{H,T}^\ell &\coloneqq \{q \in X_H \with \operatorname{supp}(q) \subset \Nb^\ell(T)\},\\
	M_{H,T}^\ell &\coloneqq \{\bbmu = (\mu, \bmu) \in M_H \with \operatorname{supp}(\mu) \subset \Sigma_T^\ell,\; \operatorname{supp}(\bmu) \subset \Nb^\ell(T)\},
\end{align*}
where \(\Sigma_T^\ell\) denotes the part of \(\Sigma\) contained in the interior of \(\Nb^\ell(T)\).
For any \(T \in \mathcal{T}_H\), a localized version of the operator \(\mathcal{K}_T\) can be defined as the map \(\mathcal{K}_T^{\ell} \colon V \rightarrow V_T^\ell\), which assigns to each \(v \in V\) the value \(\mathcal{K}_T^{\ell} v\) given by the unique solution to the following problem: find \((\mathcal{K}_T^{\ell} v, \xi_T^{\ell}, \bblambda_T^{\ell}) \in V_T^\ell \times X_{H,T}^\ell \times M_{H,T}^\ell\) such that
\begin{subequations}\label{eq:KTell}
	\begin{align}
		&\quad  a(\mathcal{K}_T^\ell v, w)& +\quad  &b(w,\xi_T^\ell)          & +\quad  &c(w,\bblambda_T^\ell) & =\quad  &-a_T(\mathcal{I}_H v, w),\quad  &\label{eq:KTella}\\
		& \quad  b(\mathcal{K}_T^\ell v,\chi)                &   &&   &             & =\quad  &-b_T(\mathcal{I}_H v,\chi),\quad &\label{eq:KTellb}\\
		&\quad c(\mathcal{K}_T^\ell v,\bbmu)  &   &                  &   &             & =\quad  &c_T(v-\mathcal{I}_H v,\bbmu)\quad&\label{eq:KTellc}
	\end{align}
\end{subequations}
for all $(w, \chi, \bbmu) \in V_T^\ell\times X_{H,T}^\ell \times
M_{H,T}^\ell$.

We now define a localized counterpart of the operator \(\mathcal{K}\) as the sum of the localized element contributions \(\mathcal{K}_T^\ell\) over all elements \(T \in \mathcal{T}_H\), cf.~\cref{eq:Ksum}:
\begin{equation}\label{eq:Kellsum}
	\mathcal{K}^{\ell} \coloneqq \sum_{T \in \mathcal{T}_H} \mathcal{K}_T^{\ell}.
\end{equation}
With this, a localized version of the \(a\)-orthogonal projection operator \(\mathcal{R}\) can be defined, recalling its decomposition in \cref{defK}, as
\begin{equation}\label{eq:Rell}
	\mathcal{R}^{\ell} \coloneqq \mathcal{I}_H + \mathcal{K}^{\ell}.
\end{equation}

The following theorem shows that the operator \(\mathcal{R}^\ell\) approximates \(\mathcal{R}\) exponentially well in the operator norm as the localization parameter \(\ell\) is increased.

\begin{lemma}[Exponential approximation]
	\label{lem:RRl}
	For all \(v \in V\) and \(\ell \in \mathbb{N}\), there holds
	\begin{equation}
		\label{eq:locerrR}
		\|\nabla (\mathcal{R} - \mathcal{R}^{\ell}) v\|_\Omega \lesssim
		\ell^{(d - 1)/2} \, \exp(- \mathfrak{c} \ell) \, \|\nabla \mathcal{R} v\|_\Omega,
	\end{equation}
	where $\mathfrak c$ is the constant from \cref{lem:expdecay}.
\end{lemma}

\begin{proof}
Let \( e \coloneqq (\mathcal{R} - \mathcal{R}^{\ell}) v \). It follows from
\cref{eq:KTc,eq:KTellc} that \( e \in W \), which, together with the definition of the
space \( Z \) in \cref{eq:defZ}, implies that \(\div e = 0\).
Thus, by the coercivity of \(a\), cf.~\cref{eq:a}, the definition of \(\tilde Z_H\) in \cref{eq:Zms}, and the  definitions of $\mathcal K^\ell$ and $\mathcal R^\ell$ in \cref{eq:Kellsum,eq:Rell}, respectively, we obtain
  \begin{equation}
    \label{eq:spliterror}  
    c_a\|\nabla (\mathcal{R}-\mathcal{R}^{\ell}) v\|_\Omega^2 \leq - a
    (\mathcal{R}^{\ell} v, e) = - \sum_{T \in \mathcal{T}_H} \big(a_T
    (\mathcal{I}_H v, e) + a (\mathcal{K}_T^{\ell} v, e)\big) .
  \end{equation}
  
Next, we estimate each term on the right-hand side of \cref{eq:spliterror} separately.  
For a fixed element \(T \in \mathcal{T}_H\), we use the cut-off function \(\eta_T\) defined as in the proof of \cref{lem:expdecay}, now with subscript \(T\).  
Noting that $a_T(\mathcal I_H v,\eta_T e) = 0$ and using the function \((1 - \eta_T)e \in V_T^{\ell}\) as a test function in \cref{eq:KTella}, yields
  \begin{align*}
    - a_T (\mathcal{I}_H v, e) - a (\mathcal{K}_T^{\ell} v, e) & = - a_T
    (\mathcal{I}_H v, (1 - \eta_T) e) - a (\mathcal{K}_T^{\ell} v, (1 -
    \eta_T) e + \eta_T e)\\
    & = - \underbrace{a (\mathcal{K}_T^{\ell} v, \eta_T e)}_{\eqqcolon \Xi_1} +
    \underbrace{b ((1 - \eta_T) e, \xi_T^{\ell})}_{\eqqcolon \Xi_2} + \underbrace{c ((1 - \eta_T) e,
    \bblambda_T^{\ell})}_{\eqqcolon \Xi_3}.
  \end{align*}
  
 For the estimate of the term \(\Xi_1\), we note that 
 \[
 \operatorname{supp}(\mathcal{K}_T^\ell v) \cap \operatorname{supp}(\eta_T e) 
 \subset R_T \coloneqq \overline{\Nb^\ell(T) \setminus \Nb^{\ell-1}(T)}.
 \]
 Recalling the \(L^\infty\)- and Lipschitz bounds of \(\eta_T\), and applying, for \(\ell \ge 2\), the Poincaré-type inequality from \cref{lem:poincare} to the functions \(e\) and \(\mathcal K_T^\ell v\) locally on the ring \(R_T\) (noting from \cref{eq:KTellc} that \(\int_E \mathcal{K}_T^\ell v \cdot n\ds = 0\) for all faces \(E \subset R_T\)), we obtain
  \begin{align*}
  	 \Xi_1  \lesssim \|\nabla \mathcal{K}_T^{\ell} v\|_{R_T} \| \nabla (\eta_T e)\|_{R_T}  + \|\mathcal K_T^\ell v\|_{R_T}\|\eta_T e\|_{R_T}\lesssim \|\nabla \mathcal{K}_T^{\ell} v\|_{R_T} \| \nabla e\|_{R_T}.
  \end{align*}
To estimate the term $\|\nabla \mathcal{K}_T^{\ell} v\|_{R_T}$ on the right-hand side of the above inequality, we apply \cref{lem:expdecay} to $\mathcal{K}_T^{\ell} v$ on $N^{\ell}(T)$ instead of $\Omega$, yielding
\begin{align}
	\label{eq:esttobecont1}
\|\nabla \mathcal{K}_T^{\ell} v\|_{R_T} \lesssim \exp(- \mathfrak{c} \ell)	\| \nabla \mathcal{K}_T^\ell v\|_{\Nb^\ell(T)}.
\end{align}
To continue the above estimate, we apply standard inf--sup theory to the saddle-point problem \cref{eq:KTell}, cf.~\cite[Cor.~4.2.1]{BoffiBrezziFortin2013}, recalling the combined inf--sup condition in \cref{infsupcbb}, now restricted to the spaces on $\Nb^\ell(T)$. This yields
\begin{multline}\label{eq:esttobecont2}
		\|\nabla K_T ^\ell v\|_{\Nb^\ell(T)} \lesssim
        \sup_{w\in V_T^\ell}\frac{|a_T(\mathcal{I}_H v, w)|}{\|\nabla w\|_\Omega}
   + \sup_{\chi\in X_{H,T}^\ell} \frac{|b_T(\mathcal{I}_H v,\chi)|}{\|\chi\|_\Omega} \\
 + \frac{1}{H}\sup _{\bbmu\in M_{H,T}^\ell} \frac{|c_T(v-\mathcal{I}_H v,\bbmu)|}{\|\bbmu\|_{M_H}}
        \\
        \lesssim  \|\mathcal I_H v\|_{1,T} + H^{-1} \|v
		-\mathcal{I}_H v\|_T + H^{-1/2} \|v -\mathcal{I}_H
		v\|_{\partial T}
		\lesssim \| v\|_{1,\Nb(T)},
\end{multline}
where we used the interpolation estimates from \cref{intEst} together with a standard trace inequality, 
and $\|\cdot\|_{1,T}$ denotes the $H^1(T)$-norm.
Noting that $\mathcal{K}_T^\ell \mathcal{R}v = \mathcal{K}_T^\ell v$, 
since $\mathcal{K}_T^\ell$ depends on $v$ only through its QOIs, which remain unchanged under the action of $\mathcal{R}$, 
and combining the previous estimates, we obtain the following bound for $\Xi_1$:
     \begin{equation}
     	\label{eq:estxi1}
     	\Xi_1 \lesssim \exp(- \mathfrak{c} \ell)	\| \mathcal R v\|_{1,\Nb(T)}\|\nabla e\|_{R_T}.
     \end{equation}
     
   To estimate the term \(\Xi_2\), we proceed analogously to the proof of \cref{lem:expdecay}, following the treatment of the corresponding term \(\Xi_2\) therein. This yields
   \[
   \Xi_2 \le \|\div\big((1-\eta_T)e\big)\|_{R}\,\|\xi_T^\ell\|_{R} 
   \lesssim \|\nabla \mathcal K_T^\ell v\|_{R_T}\,\|\nabla e\|_{R_T}.
   \]
   After applying \cref{eq:esttobecont1,eq:esttobecont2}, we obtain for \(\Xi_2\) an estimate of exactly the same form as \cref{eq:estxi1} for the term \(\Xi_1\).
Also for \(\Xi_3\), analogous arguments can be employed, constructing suitable bubble functions using \cref{LemBubbles}. The resulting estimate for \(\Xi_3\) is again of the same form as \cref{eq:estxi1}; see also the estimate of the structurally similar term \(\Xi_3\) in the proof of \cref{lem:expdecay}.

Returning to \eqref{eq:spliterror} and using the bounds for \(\Xi_1\)--\(\Xi_3\) derived above, yields
  \begin{align*}
  	\|\nabla (\mathcal{R}-\mathcal{R}^{\ell}) v\|_\Omega^2 &\lesssim \exp (- \mathfrak{c} \ell)
  	\sum_{T \in \mathcal{T}_H}  \|
  	\mathcal{R}v\|_{1,\Nb(T)} \| \nabla  	e\|_{R_T}\\
    & \lesssim \exp (- \mathfrak{c} \ell) \sqrt{\sum_{T \in \TH} \|
    \mathcal{R}v\|_{1,\Nb(T)}^2}  \sqrt{\sum_{T \in \TH} \| \nabla
    e\|_{R_T}^2}\\
    & \lesssim \ell^{(d - 1) / 2} \exp (- \mathfrak{c} \ell)  \| \nabla
    \mathcal{R}v\|_{\Omega}  \|\nabla e\|_\Omega,
  \end{align*}
where we have used the fact that each element \(K \in \mathcal{T}_H\) belongs to at most \(\mathcal{O}(\ell^{d-1})\) rings~\(R_T\) for different \(T \in \mathcal{T}_H\), as well as the Poincaré--Friedrichs inequality on~\(\Omega\). The assertion follows after dividing by \(\|\nabla e\|_{\Omega} = \|\nabla (\mathcal{R} - \mathcal{R}^{\ell}) v\|_{\Omega}\).
\end{proof}

\section{Localized multiscale method}
\label{sec:locmethod}

In this section, we introduce the proposed multiscale method for heterogeneous Stokes problems. The localized multiscale space is defined as $\tilde Z_H^\ell = \mathcal{R}^{\ell} Z$. Noting that the operator $\mathcal{R}^\ell$ depends on its argument only through the QOIs introduced in \cref{eq:qoiface,eq:qoielement}, the approximation space can be written as
\begin{equation} 
	\label{eq:approxspace} 
	\tilde Z_H^\ell \coloneqq \operatorname{span} \big\{ 
	\tilde \varphi_{F,j}^\ell \,:\, F \in \mathcal F_H^i, \, j = 1,\dots,J; \;
	\tilde \varphi_{T,k}^\ell \,:\, T \in \mathcal T_H, \, k = 1,\dots,K 
	\big\},
\end{equation}
where \(\tilde \varphi_{F,j}^\ell\) and \(\tilde \varphi_{T,k}^\ell\) are basis functions of \(\tilde Z_H^\ell\) defined as
\begin{equation}\label{eq:notusefuldefoflocbasisfun}
	\tilde \varphi_{F,j}^\ell = \mathcal{R}^\ell \tilde \varphi_{F,j},\qquad \tilde \varphi_{T,k}^\ell = \mathcal{R}^\ell \tilde \varphi_{T,k}.
\end{equation}
Although this definition involves the global prototypical basis functions, these are not needed for the actual computation of the localized basis functions, recalling that $\mathcal{R}^\ell$ depends on its input only through its readily available QOIs. A practical implementation of the proposed method, including a fine-scale discretization of the local but still infinite-dimensional patch problems in \cref{eq:KTell}, is discussed in \cref{sec:finescaledisc}.

The proposed multiscale method seeks $(\tilde u_H^\ell,\tilde p_H^\ell) \in \tilde Z_H^\ell\times Q_H$ such that
\begin{subequations}
		\label{LODid2loc} 
		\begin{align}
			\qquad \qquad \qquad & a (\tilde{u}_H^{\ell}, \tilde{v}_H^{\ell}) & + \quad
			& b (\tilde{v}_H^{\ell}, \tilde{p}_H^{\ell}) & = \quad & (f,
			\tilde{v}_H^{\ell})_{\Omega}, & \qquad \qquad \qquad\\
			\qquad \qquad \qquad & b (\tilde{u}_H^{\ell}, \tilde{q}_H^{\ell}) &  &  & =
			\quad & 0 & \qquad \qquad \qquad
		\end{align}
\end{subequations}
 for all $(\tilde{v}_H^{\ell}, \tilde{q}_H^{\ell}) \in \tilde{Z}_H^{\ell}
\times Q_H$.

The following theorem proves the well-posedness of the proposed multiscale method and its uniform convergence properties for the velocity approximations under minimal regularity assumptions, provided that the $\ell$ is chosen sufficiently large. In addition, it is proved that the piecewise constant pressure approximation~$\tilde p_H^\ell$ converges exponentially to $\Pi_H^0 p$ as the localization parameter $\ell$ is increased.

\begin{theorem}[Localized method] 
	\label{thm:convprac}
  	The localized multiscale method \cref{LODid2loc} is well-posed. Moreover, for any \( f \in H^{m+1}(\Omega) \), we have the following error estimates:
  \begin{align}
    \| \nabla (u - \tilde{u}_H^{\ell})\|_{\Omega} & \lesssim H^{m+2} |f|_{m+1, \Omega} + \ell^{(d - 1) / 2} \exp (- \mathfrak{c} \ell) \|f\|_{\Omega},\label{eq:H1errest}\\
    \|u - \tilde{u}_H^{\ell} \|_{\Omega} & \lesssim (H + \ell^{(d - 1) / 2}
    \exp (- \mathfrak{c} \ell))  \| \nabla (u - \tilde{u}_H^{\ell})\|_{\Omega},\label{eq:L2errest}\\
    \| \Pi_H^0 p - \tilde{p}_H^{\ell} \|_{\Omega} & \lesssim \ell^{(d - 1) / 2}
    \exp (- \mathfrak{c} \ell) \|f\|_{\Omega},\label{eq:pL2errest}
  \end{align}
  	where $\mathfrak c$ is the constant from \cref{lem:expdecay}.
\end{theorem}
\begin{proof}
The proofs of estimates \cref{eq:H1errest,eq:L2errest} follow exactly the same steps as in \cite[Thm.~5.1]{Hauck2025}, with the main difference being the use of the high-order error estimates from \cref{thm:convergenceprot} for the prototypical approximation. For brevity, details are omitted, and we refer to the corresponding proof in \cite{Hauck2025}. The same holds true for the proof of the pressure estimate \cref{eq:pL2errest}.
\end{proof}

\begin{remark}[Stabilized approximation]\label{rem:stabilizedapprox}
Compared to the corresponding result for the lowest-order method in \cite[Thm.~5.1]{Hauck2025}, the exponentially decaying terms in \cref{eq:H1errest}--\cref{eq:pL2errest}, representing the localization error, do not include a prefactor of \(H^{-1}\). As observed, for example, in the numerical experiments in \cite{Hauck2025}, such a prefactor can significantly degrade the quality of the localized approximation. Here, we avoid this prefactor by adapting the strategy from \cite{Hauck2025b} to the present setting. 
Note that alternative strategies exist to avoid this prefactor, as detailed, for example, in \cite{Hauck2022,Dong2023,HMM23}, and they yield quantitatively similar results in practice. However, these approaches require the bubble functions from \cref{LemBubbles} not only as a theoretical tool but also in the actual implementation. This is undesirable, especially when the bubbles are difficult to construct, as for the Stokes problem.
\end{remark}

As shown in \cref{thm:convprac}, the piecewise constant pressure approximation $\tilde p_H^\ell$ closely matches the average of the exact pressure $p$ on each mesh element, but does not capture the fine-scale oscillations present in heterogeneous settings. To address this, we introduce a post-processing step for the pressure approximation. We start by introducing  the operator \(\mathcal{R}_p^\ell\) that gathers fine-scale information about the pressure. For an input \(v \in V\), this operator is defined by
\begin{equation}\label{defRpell}
\mathcal{R}_p^\ell v \coloneqq \sum_{T \in \mathcal{T}_H} \xi_T^\ell,
\end{equation}
where \(\xi_T^\ell \in X_{H,T}^\ell\) is the second component of the solution to \cref{eq:KTell}.
The post-processed pressure approximation is defined by
\begin{equation}\label{eq:ppp}
	\tilde{p}_H^{\ell, \mathrm{pp}} \coloneqq 
	\tilde{p}_H^{\ell} + \tilde{p}_H^{\ell,\mathrm{osc}} + p_H^{\mathrm{loc}},
\end{equation}
where $\tilde{p}_H^{\ell,\mathrm{osc}} \coloneqq \mathcal{R}_p^\ell \tilde{u}_H^\ell$, and $p_H^{\mathrm{loc}} \in \mathbb{P}^{m+1}(\TH)$ is defined on each $T \in \mathcal{T}_H$ by
\begin{equation}
	\label{eq:decsourcetermproj}
	(\Pi_H^m f)|_T = \nabla p_H^{\mathrm{loc}}|_T + q_T,
\end{equation}
with $q_T \in \mathbb{Q}^m(T)$ (cf.~\cref{eq:decomposition}) and $\int_T p_H^{\mathrm{loc}} \, \mathrm{d}x = 0$.
A high-order approximation result for the post-processed pressure approximation is stated in the following theorem.
\begin{theorem}[Post-processed pressure]\label{thm:ppp}
	For any $f \in H^{m+1}(\Omega)$, we have
	\begin{equation}\label{eq:ppL2errest}
		\|p - \tilde{p}_H^{\ell, \mathrm{pp}} \|_{\Omega} \lesssim H^{m+2} |f|_{m+1,\Omega} + \ell^{(d-1)/2} \exp(-\mathfrak{c}\ell) \|f\|_{\Omega},
	\end{equation}
    	where $\mathfrak c$ is the constant from \cref{lem:expdecay}.
\end{theorem}
\begin{proof}
Consider a fixed element \(K \in \mathcal{T}_H\). Testing \cref{eq:KTella} with any \(w \in (H^1_0(K))^n\) and summing the resulting equations over all elements \(T \in \mathcal{T}_H\), while recalling that 
\(\tilde{u}_H^{\ell} = \mathcal{R}^{\ell} \tilde{u}_H^{\ell}\) and \(\tilde{p}_H^{\ell, \mathrm{osc}} = \mathcal{R}^{\ell}_p \tilde{u}_H^{\ell}\), 
yields, for some \(q_K^{\ell} \in \mathbb{Q}^m(K)\), the identity
\begin{equation}
	\label{aK1} 
	a_K(\tilde{u}_H^{\ell}, w) +b_K(w,\tilde{p}_H^{\ell, \mathrm{osc}}) = (w, q_K^{\ell})_K.
\end{equation}
Next, we use the same test function \(w\) in the reformulation of the Stokes problem~\cref{eq:reformulation}.  
Using \cref{eq:decsourcetermproj} to expand the right-hand side and noting that for $p_H = \Pi_H^0 p$ it holds $b_K(w, p_H) = -p_H|_K \int_K \div w \,\mathrm{d}x = 0$ by the divergence theorem, we obtain
\begin{equation}
	\label{aK2} 
	a_K(u, w) + b_K(w, p - {p}_H) = (f - \Pi_H^m f, w)_K + (\nabla p_H^{\mathrm{loc}} + q_K, w)_K.
\end{equation}

Taking the difference of \cref{aK1,aK2} and rearranging the terms yields
\begin{equation*}
	b_K\big(w, p - {p}_H - p_H^{\ell, \mathrm{osc}} - p_H^{\mathrm{loc}}\big) 
	= a_K(\tilde{u}_H^{\ell} - u, w) - (\tilde{q}_K^{\ell} - q_K, w)_K + (f - \Pi_H^m f, w)_K.
\end{equation*}
Applying \cref{lem:locladymod} on the element $K$ to the function
$p^\mathrm{diff} \coloneqq p - p_H - \tilde{p}_H^{\ell, \mathrm{osc}} - p_H^{\mathrm{loc}}$, which has zero mean on $K$, noting that
\begin{equation*}
	\int_K (p - {p}_H) \, \mathrm{d}x = \int_K \tilde{p}_H^{\ell, \mathrm{osc}} \, \mathrm{d}x = \int_K p_H^{\mathrm{loc}} \, \mathrm{d}x = 0,
\end{equation*}
we obtain that there exists a function \(w \in (H^1_0(K))^n\) such that 
\(-\div w = p^\mathrm{diff}\) and $(\tilde{q}_K^{\ell} - q_K, w)_K=0$.  
Using the corresponding stability estimate \(\|\nabla w\|_K \lesssim \|p^\mathrm{diff}\|_K\) together with the Poincaré–Friedrichs inequality \(\|w\|_K \lesssim H\|\nabla w\|_K\), we obtain
\begin{equation*}
	\|p^\mathrm{diff}\|_K
	\lesssim H \|\tilde{u}_H^{\ell} - u\|_K + \|\nabla (\tilde{u}_H^{\ell} - u)\|_K 
    + H \|f - \Pi_H^m f\|_K \,.
\end{equation*}
Summing the latter inequality  over all mesh elements $K\in\TH$, and applying the convergence results for the velocity approximation from \cref{thm:convprac} together with standard approximation estimates for the \(L^2\)-projection, yields
$$
\|p^\mathrm{diff}\|_\Omega
\lesssim H^{m+2} |f|_{m+1, \Omega} + \ell^{(d - 1) / 2} \exp (- \mathfrak{c} \ell) \|f\|_{\Omega} \,.
$$
To derive the desired estimate \cref{eq:ppL2errest}, we write $p - \tilde{p}_H^{\ell, \mathrm{pp}} = \Pi_H^0 p - \tilde{p}_H^\ell + p^\mathrm{diff}$ and apply the triangle inequality, combining the above estimate with \cref{eq:pL2errest}.
\end{proof}

\section{Fine-scale discretization}\label{sec:finescaledisc}

In this section, we discuss the practical implementation of the proposed multiscale method, including the computation of the localized basis functions via a fine-scale discretization of the local infinite-dimensional patch problems in \cref{eq:KTellprac}.
Let \(\mathcal{T}_h\) be a mesh of \(\Omega\) fine enough to resolve all microscopic features of the coefficients, and let \(V_h \subset V\) denote the corresponding finite element space for the velocity.
 We assume that the fine mesh \(\mathcal{T}_h\) is compatible with the coarse mesh \(\mathcal{T}_H\), in the sense that the restriction of \(\mathcal{T}_h\) to any element \(T \in \mathcal{T}_H\), denoted \(\mathcal{T}_h(T)\), forms a valid mesh itself.
For any \(T \in \mathcal{T}_H\), let~\(V_h(T)\) denote the restriction of \(V_h\) to \(T\), that is, $V_h(T) \coloneqq \{ v_h|_T : v_h \in V_h\},$ and let \(Q_h(T) \subset L^2(T)\) be the associated local pressure space.  
The global pressure approximation space is then defined by
\[
Q_h = \bigl\{ q_h \in Q \with  q_h|_T \in Q_h(T),\; \forall T \in \mathcal{T}_H \bigr\}.
\]

We assume that the discrete spaces satisfy the following natural conditions. First, the discrete velocity–pressure inf--sup condition should hold not only globally between the spaces \(V_h\) and \(Q_h\), but also locally on each coarse element \(T\) between 
\(V_h^0(T) = V_h(T) \cap H^1_0(T)^n\) and \(X_h(T) = Q_h(T) \cap L^2_0(T)\), that is,
\begin{equation} \label{infsuploc} 
\adjustlimits \inf_{q_h \in X_h(T)} \sup_{v_h \in V_h^0(T)} 
\frac{b_T(v_h, q_h)}{\|\nabla v_h\|_T\|q_h\|_T } \gtrsim 1.  
\end{equation}       
Second, the spaces provide order \(k+1\) approximation for velocity and pressure:
\begin{align*}
\forall v \in H^1_0(\Omega) \cap H^{k+2}(\Omega), \quad
& \inf_{v_h \in V_h} \big(\|\nabla(v - v_h)\|_\Omega + h^{-1} \|v - v_h\|_\Omega \big)
\lesssim h^{k+1} |v|_{k+2,\Omega}, \\
\forall q \in L^2_0(\Omega) \cap H^{k+1}(\Omega), \quad
& \inf_{q_h \in Q_h} \|q - q_h\|_\Omega \lesssim h^{k+1} |q|_{k+1,\Omega}.
\end{align*}

\begin{samepage}
    This convention covers the two cases:
\begin{itemize}
	\item Velocity–pressure pairs with a discontinuous pressure approximation, such as the Scott--Vogelius pair \(\mathcal{P}^{k+1}_\mathrm{cg}/\mathcal{P}_{\mathrm{dg}}^k\), cf.~\cite{Scott1985}. For the barycentric refinement of a mesh, both the global inf--sup stability and its local version~\cref{infsuploc}  hold for any polynomial degree \(k \ge 1\); see, e.g.,~\cite{Guzman2018}.
	
	\item Velocity–pressure pairs with pressure continuous on each coarse element, closely related to the Taylor--Hood pair $\mathcal P_\mathrm{cg}^{k+1}/\mathcal P_\mathrm{cg}^k$; cf.~\cite{Girault1986}.  
The global pair \(V_h / Q_h\) allows discontinuities of the pressure across coarse elements and is therefore not a standard Taylor--Hood pair.  
Nevertheless, both the global inf--sup condition for \(V_h / Q_h\) and its local version~\cref{infsuploc} follow from the stability of the local pairs \(V_h(T) / Q_h(T)\), which is guaranteed by standard Taylor--Hood theory; we refer to~\cref{sec:InfSupTHmod} for a proof.
\end{itemize}
\end{samepage}

In both cases, the space $Z_h$, a fully discrete counterpart of $Z$ from \cref{eq:defZ} appearing in the reformulation \cref{eq:reformulation} of the Stokes problem, can be defined as
\[
Z_h = \{ v_h \in V_h \with b(v_h, q_h) = 0 \ \forall q_h \in X_{H,h} \},
\] 
where the fine-scale discretized version of $X_H$ is given by
\[
X_{H,h} = \{ q_h \in Q_h \with \Pi_H^0 q_h = 0 \}.
\]
With these spaces at hand, we can define fully discrete versions of the fine-scale 
space $W_h$ and the prototypical approximation space $\tilde Z_{H,h}$, 
analogous to \cref{eq:defW,eq:Zms}, respectively. Moreover, a fine-scale 
discretization of the prototypical method can be defined analogous to 
\cref{LODid2}.

To define a fully discrete version of the localized method \cref{LODid2loc}, the local patch problems for the element contributions $\mathcal K_T^\ell$ in \cref{eq:KTell} must be discretized. For this, we introduce fully discrete counterparts of the spaces $V_T^\ell$ and $X_{H;T}^\ell$ in \cref{eq:KTell}, namely
\[
V_{T,h}^\ell = V_T^\ell \cap V_h, \qquad 
X_{H,T,h}^\ell \coloneqq X_{H,T}^\ell \cap X_{H,h}.
\]
Having defined these spaces, a fully discrete version of the operator $\mathcal R^\ell$ from~\cref{eq:Rell} can be introduced, and the fully discrete basis functions can be constructed analogously to \cref{eq:notusefuldefoflocbasisfun}. To make this definition suitable for a practical use, we reformulate them taking into the account that the fully discrete localized basis functions should have the same QOIs as the corresponding prototypical basis functions. 

For example, the fully discrete basis function $\tilde \varphi^\ell_{F,j,h}$ associated with a face $F$ and an index $j\in\{1,\ldots,J\}$ should satisfy, for all $\bbmu \in M_H$,
 $$c(\tilde \varphi^\ell_{F,j,h},\bbmu)=H\int_F p_{F,j}\mu\ds.$$
To compute the fully discrete basis function $\tilde \varphi^\ell_{F,j}$ in practice, we first introduce coefficients $\theta_z \in \mathbb{R}^n$ for each interior node $z$ such that
\[
\mathcal{I}_H \tilde \varphi^\ell_{F,j,h} = \sum_{z} \theta_z \Lambda_z,
\] 
where $\Lambda_z$ denotes the standard (scalar-valued) finite element hat function associated with the node $z$. The coefficients $\theta_z$ can be directly obtained from the definition of the operator $\mathcal{I}_H$, which depends only on the normal face integrals of its argument. Indeed, for any face $E \in \mathcal{F}_H^i$, we have $|E|^{-1}\int_E \tilde \varphi^\ell_{F,j,h} \cdot n \ds = \delta_{EF}\delta_{j1}$. Although the coefficients $\theta_z$ depend on $F$ and $j$, this dependence is omitted for notational simplicity (and likewise for the forthcoming functions $\psi_{T,h}^\ell$).
The fully discrete localized basis function, together with the corresponding pressure contribution (to be used in the implementation of operator $\mathcal{R}_p^\ell$ in (\ref{defRpell}), needed for the pressure post-processing), can then be computed as
\begin{equation} 
	\label{eq:locbasis}
	\tilde{\varphi}^\ell_{F,j,h} = \sum_z \theta_{z} \Lambda_z + \sum_{T\in\TH} \psi_{T,h}^\ell,  \quad
	\tilde{\xi}^\ell_{F,j,h} = \sum_{T\in\TH} \xi_{T,h}^\ell,
\end{equation}
where \( (\psi_{T,h}^\ell, \xi_{T,h}^\ell, \bblambda_{T}^\ell) \in V_{T,h}^\ell \times X_{H,T,h}^\ell\times M_{H,T}^\ell \) solves the saddle-point problem
\begin{subequations}\label{eq:KTellprac}
	\begin{align}
		&a(\psi_{T,h}^\ell, w_h)& +\quad  &b(w_h,\xi_{T,h}^\ell)          & +\quad  &c(w_h,\bblambda_{T}^\ell) & =&- \textstyle \sum_{z}  a_{T}(\theta_z\Lambda_z, w_h),\quad  &\label{eq:KTellpraca}\\
		&b(\psi_{T,h}^\ell,\chi_h)                &   &&   &             & = &- \textstyle \sum_{z}  b_{T}(\theta_z\Lambda_z,\chi_h),\quad &\label{eq:KTellpracb}\\
		&c(\psi_{T,h}^\ell,\bbmu)  &   &                  &   &             & = &\ g_{T}(\bbmu)\quad&\label{eq:KTellpracc}
	\end{align}
\end{subequations}
for all $(w_h, \chi_h, \bbmu) \in V_{T,h}^\ell\times X_{H,T,h}^\ell \times
M_{H,T}^\ell$, and
\begin{equation*}
	g_T(\bbmu) \coloneqq H   \tint_F \kappa_T p_{F,j} \mu \ds - H \textstyle \sum_{z}  \int_{\partial T} \kappa_T \Lambda_z\theta_z\cdot n \mu \ds - \sum_{z}  \int_T (\theta_z\Lambda_z) \cdot \bmu \dx
\end{equation*}
with unit normal \(n \colon \Sigma \to \mathbb{R}^n\)  defined in \cref{eq:defn}. 
We emphasize that only a few local problems of the form \cref{eq:KTellprac} need to be solved for each basis function. For $j=1$, the functions~$\psi_{T,h}^\ell$ are nonzero only for elements $T$ sharing a node with $F$. Furthermore, for $j>1$, they are nonzero only for elements $T$ sharing the face $F$.

For an element $T \in \TH$ and index $k \in \{1,\ldots,K\}$, the fully discrete basis functions $\tilde\varphi_{T,k,h}^\ell$ and their corresponding pressure contributions $\tilde\xi_{T,k,h}^\ell$ are computed similarly, but more simply, since only one problem of the form \cref{eq:KTellprac} must be solved per basis function. Since $\mathcal{I}_H \tilde\varphi_{T,k,h}^\ell = 0$, both $\tilde\varphi_{T,k,h}^\ell$ and $\tilde\xi_{T,k,h}^\ell$ can be obtained by solving \cref{eq:KTellprac} directly with $\theta_z = 0$ for all nodes $z$ and $g_T(\bbmu) = \int_T p_{T,k} \cdot \bbmu \,\mathrm{d}x$.

The fully discrete approximation space is then given by
\begin{equation*} 
	\label{heq:approxspace} 
	\tilde Z_{H,h}^\ell \coloneqq \operatorname{span} \big\{ 
	\tilde \varphi_{F,j,h}^\ell \,:\, F \in \mathcal F_H^i, \, j = 1,\dots,J; \;
	\tilde \varphi_{T,k,h}^\ell \,:\, T \in \mathcal T_H, \, k = 1,\dots,K 
	\big\},
\end{equation*}
and the fully discrete multiscale method then seeks $(\tilde u_{H,h}^\ell, \tilde p_{H,h}^\ell) \in \tilde Z_{H,h}^\ell \times Q_H$, solving \cref{LODid2loc} with the corresponding notational changes. Its solution is given by
$$
  \tilde u_{H,h}^\ell =
	\sum_{F \in \mathcal F_H^i}\sum_{j=1}^J u_{F,j}\tilde \varphi_{F,j,h}^\ell
    \,+\,
	\sum_{T \in \TH}\sum_{k=1}^Ku_{T,k}\varphi_{T,k,h}^\ell,
$$
for some coefficients $\{u_{F,j}\}_{F,j}$ and $\{u_{T,k}\}_{T,k}$. These coefficients can then be used to compute the oscillatory pressure contribution, cf. \cref{defRpell}--\cref{eq:ppp}, by setting
$$
  \tilde{p}_{H,h}^{\ell,\mathrm{osc}} =
	\sum_{F \in \mathcal F_H^i}\sum_{j=1}^J u_{F,j}\tilde \xi_{F,j,h}^\ell
    \,+\,
	\sum_{T \in \TH}\sum_{k=1}^K u_{T,k}\xi_{T,k,h}^\ell \,.
$$
The post-processed pressure is then defined as the fully discrete analogue of \cref{eq:ppp}:
$$
	\tilde{p}_{H,h}^{\ell, \mathrm{pp}} \coloneqq 
	\tilde{p}_{H,h}^{\ell} + \tilde{p}_{H,h}^{\ell,\mathrm{osc}} + p_H^{\mathrm{loc}}.
$$

For this fully discrete approximation, a convergence result analogous to \cref{thm:convprac} is stated below, with the error measured against the fine-scale finite element reference solution $(u_h, p_h)$ obtained from \cref{eq:weakstokes} using the pair $V_h/Q_h$.

\begin{theorem}[Fully discrete localized method]\label{thm:fd}
    There exists \(\gamma > 0\) depending only on the regularity of meshes $\TH$ and $\Th$, such that the fully discrete version of the localized multiscale method \cref{LODid2loc} is well-posed provided $\tfrac{h}{H} \leq \gamma$. Moreover, for any \( f \in H^{m+1}(\Omega) \), we have the following error estimates:
	\begin{align}
		\|\nabla ( u_h - \tilde{u}_{H, h}^{\ell})\|_{ \Omega} &\lesssim H^{m + 2} | f |_{m + 1,
			\Omega} + h^{k+1} | f |_{k, \Omega}  + \ell^{(d - 1) / 2} \exp (- \mathfrak{c} \ell) \|f\|_{\Omega},\label{heq:H1errest}\\
        \|u_h - \tilde{u}_{H, h}^{\ell} \|_{\Omega} & \lesssim (H + \ell^{(d - 1) / 2}
    \exp (- \mathfrak{c} \ell))  \|\nabla ( u_h - \tilde{u}_{H, h}^{\ell})\|_{ \Omega},\label{heq:L2errest}\\
			    \|p_h - \tilde{p}_{H,h}^{\ell, \mathrm{pp}} \|_{\Omega}  &\lesssim H^{m + 2}|f|_{m + 1, \Omega} + h^{k+1} | f |_{k, \Omega} +\ell^{(d - 1) / 2} \exp (- \mathfrak{c} \ell) \|f\|_{\Omega}.	\label{heq:ppL2errest}
	\end{align}
The term $h^{k+1} |f|_{k,\Omega}$ in the above error estimates can be omitted if the fine-scale discretization is divergence-free or if $k \ge m+1$.
\end{theorem}
\begin{proof}
We recall that the proof of the error estimates in \cref{thm:convprac,thm:ppp} boils down to the following essential steps: the well-posedness of the saddle-point problems for the prototypical basis functions in \cref{le:protbasis}, the error estimate for the prototypical method in \cref{thm:convergenceprot}, the proof of the exponential decay in \cref{lem:expdecay}, and the estimate of the error due to localization in \cref{lem:RRl}. Examining the proofs of these results, one realizes that they are all based on the existence of bubble functions, as constructed in \cref{lem:locladymod,LemBubbles}. In order to transfer these proofs to the fully discrete setting, one should thus first construct the discrete analogues of these Lemmas. This is done in \cref{AppendixBubble}, cf. \cref{hlem:locladymod,hLemBubbles}.  Using the fully discrete bubble functions from these Lemmas, one can directly recast the proofs of \cref{le:protbasis,lem:expdecay,lem:RRl} to the fully discrete setting, replacing the continuous functions spaces by their discrete counterparts.

It remains to adapt the proof of \cref{thm:convergenceprot}, i.e. to bound the error of the fully discrete prototypical approximation $\tilde u_{H,h}$ against the fine-scale velocity approximation $u_h$. For this, we revisit the proof of \cref{thm:convergenceprot}, setting $e \coloneqq u_h - \tilde{u}_{H,h}$. 
The only modification in the proof occurs in~\cref{trickyCR}, where we must account for $\div e \neq 0$ which is the case if the fine-scale discretization is not exactly divergence-free (otherwise, the proof carries over unchanged). Thus, \cref{trickyCR} becomes
	\[ (\Pi_H^m f, e)_{\Omega} = \sum_{F \in \mathcal{F}_H^i} ([\phi]_F, e \cdot
	n)_F - \sum_{T \in \mathcal{T}_H} (\phi_T - \phi_{T, h}, \ \div e)_T,
	\]
	Here we could subtract any $\phi_{T,h} \in Q_h(T)$ since $b(\tilde{u}_{H,h},q_h) = 0$ for all $q_h \in Q_h$. Note that $\phi_T$ is a polynomial of degree $m+1$, while $\phi_{T,h}$ is a piecewise polynomial of degree $k$. Hence, the second term above vanishes if $k \ge m+1$.
    
Otherwise, proceeding as in the continuous case in the proof of \cref{thm:convergenceprot} and additionally choosing $\phi_{T,h}$ as a suitable interpolant of $\phi_T$, we conclude
    \begin{align*}
         | (\Pi_H^m f, e)_{\Omega} | &\lesssim H^{m + 2} | f |_{m + 1, \Omega} \|\nabla e
	\|_{\Omega} + h^{k+1} \sum_{T \in \mathcal{T}_H} | \phi_T |_{k+1, T} \|\nabla  e \|_{
		T}\\
         &\lesssim H^{m + 2} | f |_{m + 1, \Omega} \|\nabla e \|_{\Omega} + h^{k+1} \sum_{T\in \TH} | g_T |_{k, T} \|\nabla e \|_{T} \\
        &\lesssim  (H^{m + 2} |
	f |_{m + 1, \Omega} + h^{k+1} | f |_{k, \Omega}) \| \nabla e \|_{\Omega}
    \end{align*}
This leads to the discrete analogue of \cref{thm:convergenceprot}. Putting all the above mentioned ingredients together leads to the estimates  (\ref{heq:H1errest})--(\ref{heq:ppL2errest}).
\end{proof}

The error estimates in \cref{thm:fd} are against the fine-scale finite element solution. To obtain an error estimate against the continuous solution to problem \cref{eq:weakstokes}, we apply the triangle inequality, together with the classical convergence result
\[
\|\nabla(u - u_h)\|_{\Omega} + \|p - p_h\|_{\Omega} \lesssim h^s \big( |u|_{1+s,\Omega} + \|p\|_{s,\Omega} \big),
\]
where we make the regularity assumption $u \in H^{1+s}(\Omega)$ and $p \in H^s(\Omega)$ for some parameter $0 < s \le k+1$. We emphasize that, for heterogeneous Stokes problems, the solution seminorms on the right-hand side may be large or the regularity parameter~$s$ may be close to zero, leading to reduced convergence rates.

\section{Numerical experiments}\label{sec:numexp}

\begin{figure}
	\centering
	\begin{subfigure}[t]{0.32\textwidth}
		\centering
		\includegraphics[height=.89\linewidth]{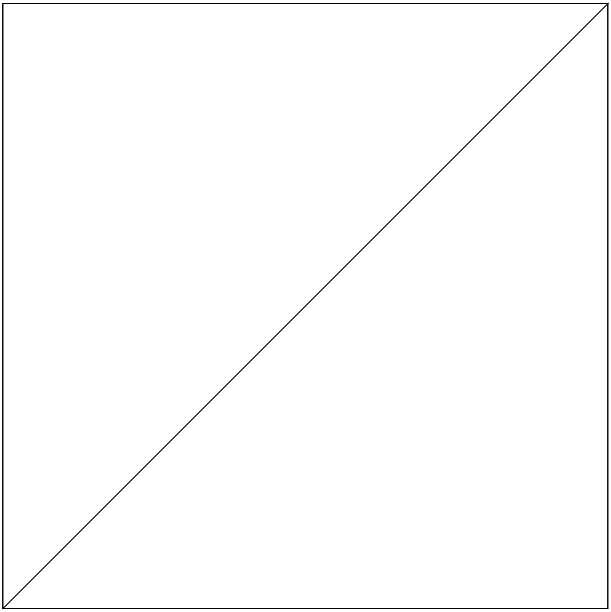}
	\end{subfigure}
	\begin{subfigure}[t]{0.32\textwidth}
		\centering
		\includegraphics[height=.89\linewidth]{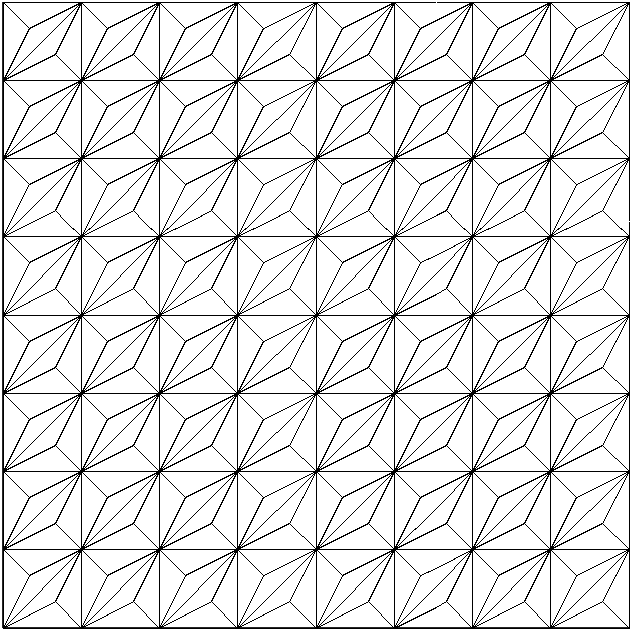}
	\end{subfigure}
	~ 
	\begin{subfigure}[t]{0.32\textwidth}
		\centering
		\includegraphics[height=.9\linewidth]{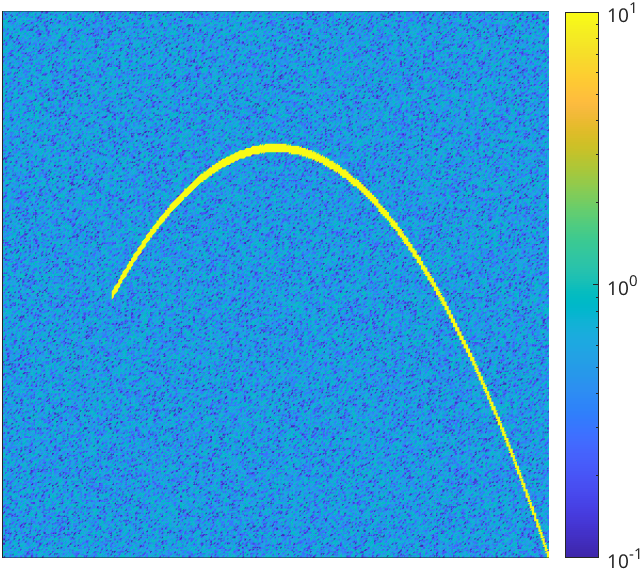}
	\end{subfigure}
	\caption{Initial mesh $\mathcal T_{2^{-0}}$ for the mesh generation (left), barycentric refinement of mesh $\mathcal T_{2^{-3}}$, and multiscale coefficient used in all numerical experiments (right).}
	\label{fig:coeff}
\end{figure}

In this section, we numerically investigate the proposed multiscale method. For all numerical experiments, we consider the domain $\Omega = (0,1)^2$ and a hierarchy of meshes generated by uniform red refinement 
of the initial mesh shown in \cref{fig:coeff}~(left). For simplicity, the meshes 
in the hierarchy are denoted by $\mathcal T_{2^0}, \mathcal T_{2^{-1}}, \dots$, 
where the subscript indicates the side length of the squares formed by joining 
opposing triangles.  
We choose the viscosity coefficient $\nu$ of the Stokes problem \cref{pbStokes} 
 to be piecewise constant on the mesh $\mathcal T_\epsilon$ 
with $\epsilon = 2^{-8}$, and we draw its element values as independent uniform 
random variables in $[0.1,1]$. For elements whose midpoints lie within a distance 
$4\epsilon$ of some prescribed parabola, the corresponding element values are set to 10; see 
\cref{fig:coeff} (right). The damping coefficient $\sigma$ is set to zero for 
simplicity. The source term is chosen as
\begin{equation*}
	f(x,y) \coloneqq (-y,x^4)^\top.
\end{equation*}

For the fine-scale discretization, we use the Scott--Vogelius pair 
$\mathcal P^{k+1}_\mathrm{cg}/\mathcal P^k_\mathrm{dg}$ with $k=1$ on a mesh 
obtained by uniform barycentric refinement of $\mathcal T_{2^{-8}}$. 
\cref{fig:coeff} (center) illustrates such a refinement using a rather coarse mesh for clarity. 
For barycentrically refined meshes, the inf--sup stability of the Scott--Vogelius 
element is guaranteed for any polynomial degree $k \ge 1$; see, e.g.,~\cite{Guzman2018}.
Note that all numerical experiments presented below can be reproduced using the code available at  \url{https://github.com/moimmahauck/Stokes_HO_LOD}. 

\subsection{High-order convergence}

First, we study the convergence of the proposed high-order multiscale 
method~\cref{LODid2loc} under mesh refinement. To this end, we introduce the 
following error measures for the velocity and pressure approximations:
\begin{align*}
	\mathrm{err}_{u,H^1}(H,\ell) &\coloneqq \|\nabla (u_h - \tilde u_{H,h}^\ell)\|_\Omega, 
	&\qquad 
	\mathrm{err}_{u,L^2}(H,\ell) &\coloneqq \|u_h - \tilde u_{H,h}^\ell\|_\Omega,\\[0.2em]
	\mathrm{err}_{p,L^2}(H,\ell) &\coloneqq \|p_h - \tilde p_{H,h}^{\ell,\mathrm{pp}}\|_\Omega, 
	&\qquad
	\mathrm{err}_{\Pi_H p,L^2}(H,\ell) &\coloneqq \|\Pi_H p_h - \tilde p_{H,h}^\ell\|_\Omega.
\end{align*}

\begin{figure}
	\includegraphics[height=.32\linewidth]{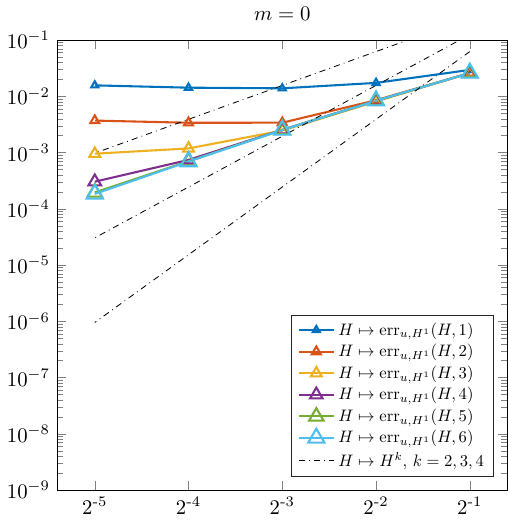}\hfill
	\includegraphics[height=.32\linewidth]{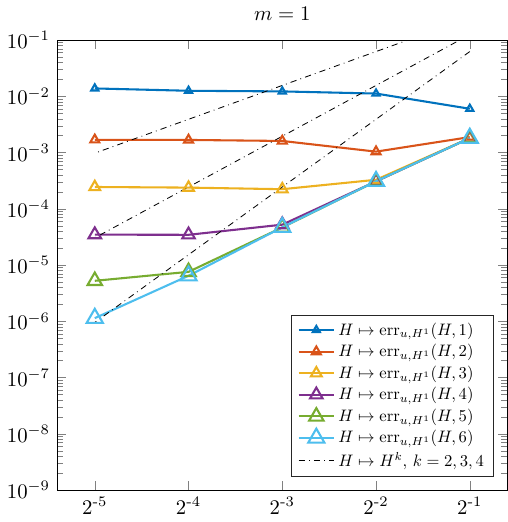}\hfill
	\includegraphics[height=.32\linewidth]{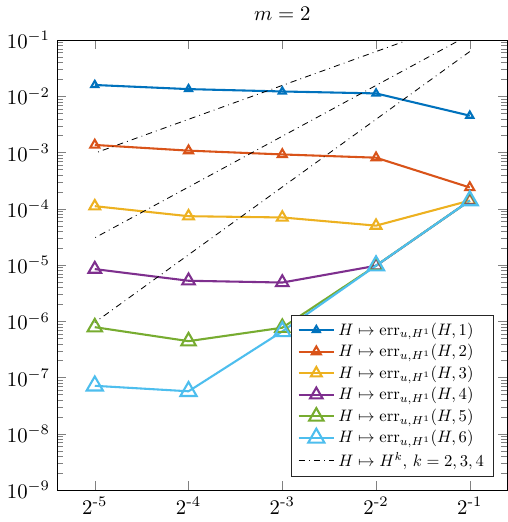}\\[2ex]
	\includegraphics[width=.32\linewidth]{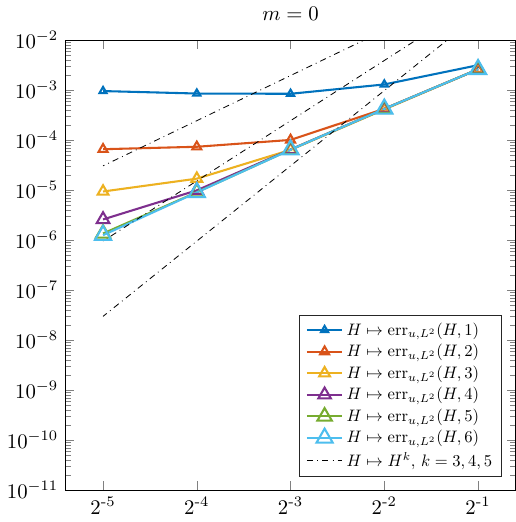}\hfill
	\includegraphics[width=.32\linewidth]{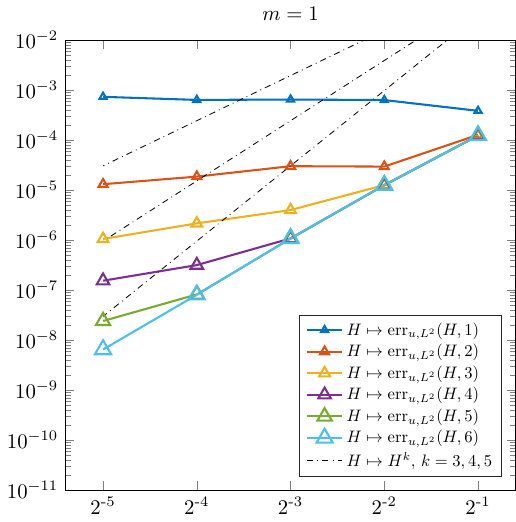}\hfill
	\includegraphics[width=.32\linewidth]{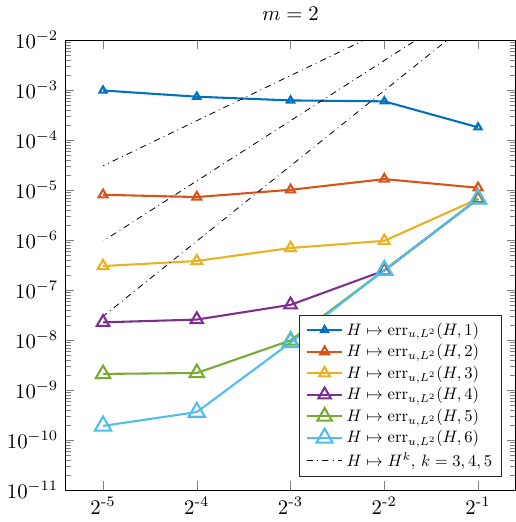}
	\caption{Error plots of the velocity approximation for polynomial degrees $m \in \{0,1,2\}$ (from left to right). For fixed localization parameters $\ell$, the $H^1$-norm (top row) and $L^2$-norm (bottom row) errors are plotted as functions of the coarse mesh size $H$.
	}
	\label{fig:convu}
\end{figure}

For the $L^2$- and $H^1$-errors of the velocity approximation, \cref{fig:convu} 
shows convergence orders of $m+2$ and $m+3$, respectively, provided the 
localization parameter is sufficiently large. Since the source term $f$ is smooth, 
these convergence rates agree with the theoretical predictions from \cref{thm:convprac}. For a fixed localization parameter, we observe that, once the error reaches a error level, 
it stagnates as the mesh is further refined. This error level is determined by the localization error for the chosen localization parameter. The distance between these plateaus 
increases for higher polynomial degrees, indicating that the localization 
properties of the method improve with increasing polynomial order. This behavior 
is consistent with earlier observations for elliptic diffusion-type problems; see, e.g., 
\cite{Maier2021,Dong2023,Hauck2025b}.

\begin{figure}
	\includegraphics[height=.32\linewidth]{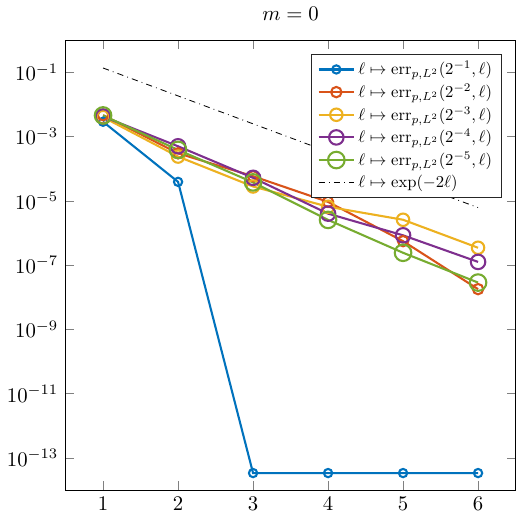}\hfill
	\includegraphics[height=.32\linewidth]{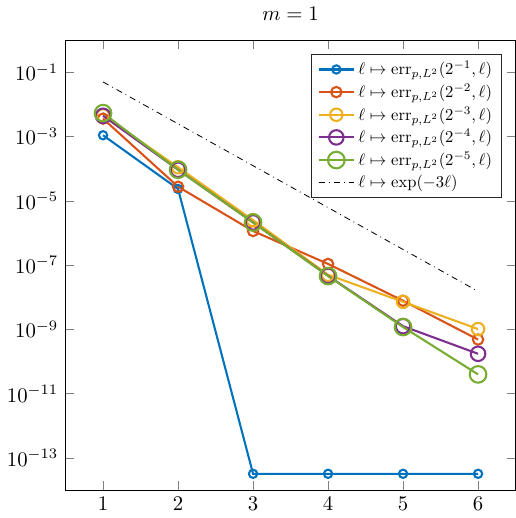}\hfill
	\includegraphics[height=.32\linewidth]{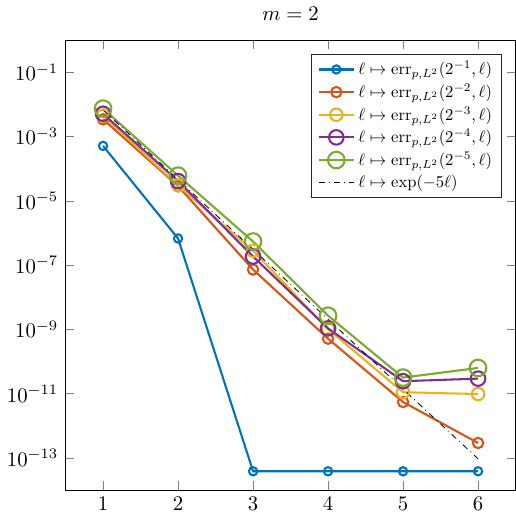}
	\caption{Error plots of the pressure approximation for polynomial degrees 
		$m \in \{0,1,2\}$ (from left to right). For fixed coarse mesh sizes $H$, the 
		$L^2$-norm error relative to $\Pi_H p_h$ is shown as a function of the 
		localization parameter $\ell$.}
	\label{fig:convp}
\end{figure}

In \cref{fig:convp}, we observe that the piecewise constant pressure approximation 
of the proposed method converges exponentially towards $\Pi_H p_h$, in agreement 
with the theoretical prediction of \cref{thm:fd}. 
Note that the blue curves reach machine accuracy for localization parameters $\ell \geq 3$. This is because, for the corresponding mesh size $H = 2^{-1}$ and $\ell \geq 3$, the patches on which the basis functions are defined already cover the entire domain $\Omega$, so the exponentially decaying localization error is zero.
In \cref{fig:convp}, one also observes the improved localization properties obtained by increasing 
the order $m$ of the method. Note that for the dashed, exponentially decaying 
reference lines, the decay rates increase as $m$ becomes larger.

\begin{figure}
		\includegraphics[height=.32\linewidth]{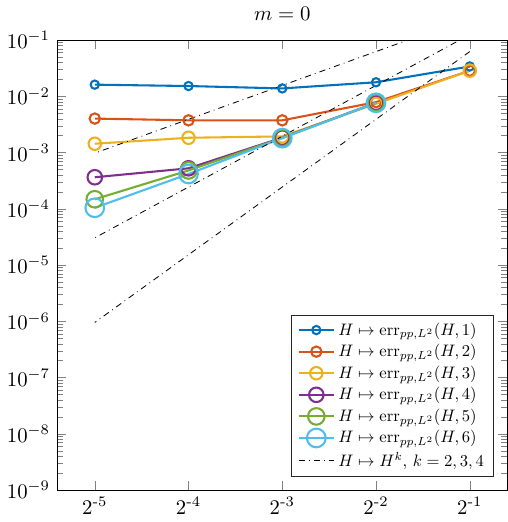}\hfill
		\includegraphics[height=.32\linewidth]{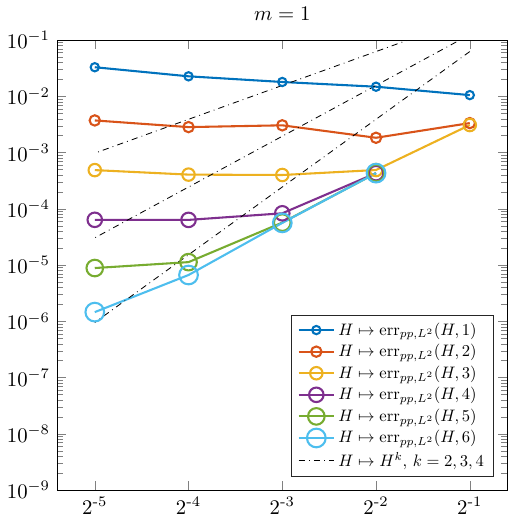}\hfill
		\includegraphics[height=.32\linewidth]{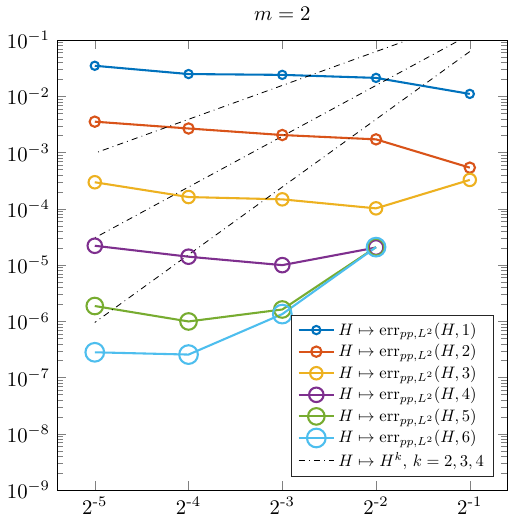}
		\caption{Error plots of the post-processed pressure approximation for polynomial degrees 
		$m \in \{0,1,2\}$ (from left to right). For fixed localization parameters $\ell$, the $L^2$-norm errors are plotted as functions of the coarse mesh size $H$.}
		\label{fig:convppp}
\end{figure}

\cref{fig:convppp} shows that the post-processing step yields a highly accurate 
pressure approximation, despite originating from an inaccurate piecewise constant coarse-scale 
 approximation. The error levels and convergence 
behavior are, for all orders~$m$, comparable to those of the $H^1$-approximation of 
the velocity, which converges with order $m+2$. This is consistent with the 
theoretical prediction of \cref{thm:fd}.

\subsection{Comparison with lowest-order method from \cite{Hauck2025}}

\begin{figure}
	\begin{center}
		\includegraphics[height=.32\linewidth]{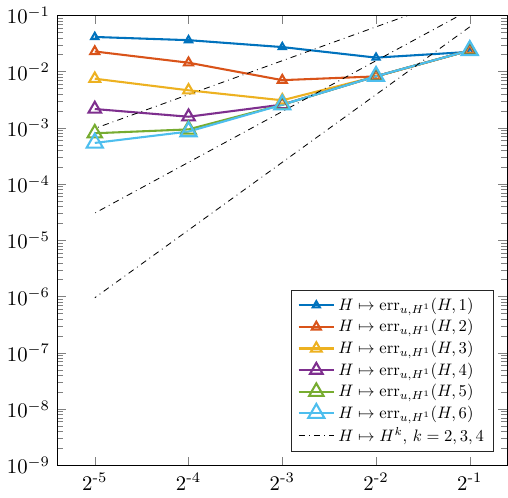}\hspace{2ex}
		\includegraphics[height=.32\linewidth]{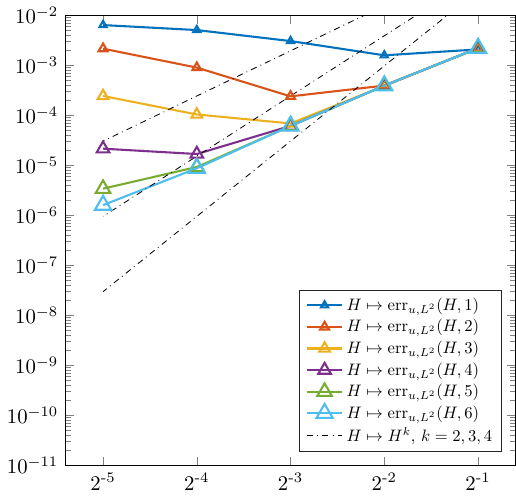}\\[2ex]
		\includegraphics[height=.32\linewidth]{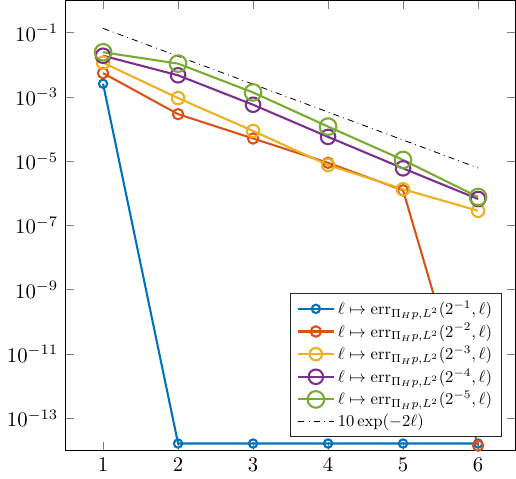}\hspace{2ex}
		\includegraphics[height=.32\linewidth]{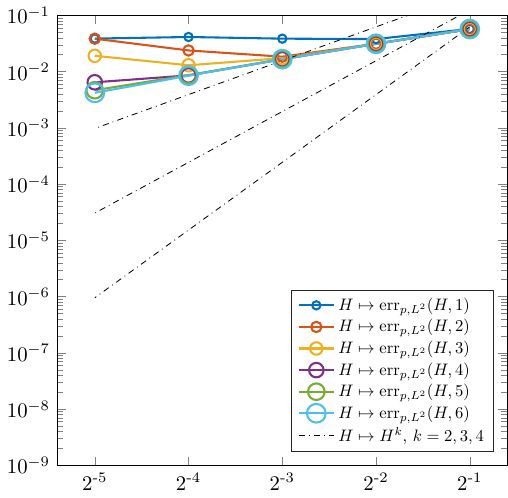}
	\end{center}
	\caption{Error plots for the lowest-order multiscale method from \cite{Hauck2025}: 
		velocity $H^1$-error, velocity $L^2$-error, pressure $L^2$-error with respect to 
		$\Pi_H p_h$, and post-processed pressure $L^2$-error (arranged from left to right, 
		top to bottom).}
	\label{fig:convcomparison}
\end{figure}

Finally, we numerically compare the proposed method with the lowest-order 
multiscale method introduced in \cite{Hauck2025} for heterogeneous Stokes 
problems. We consider exactly the same problem setup, including the domain, 
coefficient, source term, and meshes. Note, however, that the fine-scale 
discretization differs: in \cite{Hauck2025} a Crouzeix--Raviart method is used to this end. This difference does not affect the 
quantitative behavior of the observed errors, since the error of 
the multiscale method is also measured with respect to the 
fine-scale Crouzeix--Raviart~solution.

The main difference between the multiscale method in \cite{Hauck2025} and the
lowest-order version of the present method ($m = 0$) lies in the choice of
QOIs: \cite{Hauck2025} uses edge integrals, whereas our method relies solely on
normal face integrals; see \cref{eq:qoiface}. This reduces the dimension of the
approximation space by a factor of two in two dimensions. Despite this smaller
space, the achievable error levels for large localization parameters are
comparable for both lowest-order methods; cf.\ \cref{fig:convcomparison}. Their
localization properties also appear to be similar.

For fixed localization parameters, however, the method from \cite{Hauck2025} exhibits 
the issue that, after reaching a certain error level, the error increases again as the 
mesh is further refined. This phenomenon is well-known in the LOD community 
(see also \cref{rem:stabilizedapprox} and \cite{MalP14,Maier2021}) and is caused by the rather naive localization strategy 
used in these methods. In contrast, the more sophisticated localization strategy employed 
here, relying on a quasi-interpolation operator $\mathcal I_H$ (cf. \cref{intEst}) and defining the operator 
$\mathcal R^\ell$ as a sum of element contributions (cf. \cref{eq:Rell,eq:Kellsum}), eliminates this effect. For the 
proposed method, the error stagnates as the mesh is refined, once the localization 
parameter-dependent error level is reached.

Another improvement compared to the method in \cite{Hauck2025}, which also appears in the lowest-order case of the method proposed here, is the refined post-processing strategy. In \cite{Hauck2025}, post-processing consisted solely of adding the oscillatory term $\tilde p_H^{\mathrm{osc},\ell}$ (see \cref{eq:ppp}) to the piecewise constant pressure approximation. As a result, as proved in \cite[Thm.~5.1]{Hauck2025}, the post-processed pressure achieves only first-order convergence; see \cref{fig:convcomparison} (bottom left) for a numerical confirmation.  
In this work, we prove that by adding an additional coarse, locally polynomial term, denoted by~$p^\mathrm{loc}$ in \cref{eq:ppp}, the post-processed pressure can attain the same convergence order as the $H^1$-velocity approximation (see \cref{thm:ppp}). In the lowest-order case $m = 0$, this corresponds to second-order convergence.

\section*{Acknowledgments}
M.~Hauck acknowledges funding from the Deutsche Forschungsgemeinschaft\linebreak  (DFG, German Research Foundation) -- Project-ID 258734477 -- SFB 1173. 

\appendix

\section{Proofs of some technical lemmas}\label{AppendixBubble}

\begin{proof}[Proof of \cref{lem:poincare}]
	Note that if $v \in (H^1(T))^n$ satisfies $\tint_F v\cdot n \ds  = 0$ for all faces $F \subset \partial T$, then $\|\nabla v\|_T=0$ implies that $v=0$. More precisely, such $v$ is a constant vector satisfying $v\cdot n=0$ for all normal vectors $n$ associated with  faces $F \subset \partial T$. Since we have at least $n$ linearly independent vectors for each element $T$, this implies that $v=0$. This observation together with the Peetre-Tartar lemma, cf.~\cite[Lem.~A.38]{ErG04}, and a scaling argument concludes the assertion. 
\end{proof}

We now turn to the proofs of \cref{lem:locladymod,LemBubbles}. Both results follow as corollaries of the lemma stated below.

\begin{lemma}[Compliance with element constraints]
	\label{LemCelWeights}
	For any element $T \in \mathcal{T}_H$ and any linear functional 
	$L \colon (\mathbb{P}^m(T))^n \to \mathbb{R}$ satisfying 
	$L g = 0$ for all $g \in \mathbb{G}^m(T)$, there exists a local function 
	$v_L \in (H^1_0(T))^n$ such that $\div v_L = 0$,
	\begin{equation}
		\label{defbT}
		\int_T v_L \cdot w \, \mathrm{d}x = L w, \quad \forall w \in (\mathbb{P}^m(T))^n,
	\end{equation}
	and the following stability estimate holds:
	\begin{equation}
		\label{estbT}
		\|v_L\|_T + H \|\nabla v_L\|_T \lesssim  \|L\|, \qquad 
		\|L\| \coloneqq \sup_{{w \in (\mathbb{P}^m(T))^n \with \|w\|_T = 1}} \|L w\|_T.
	\end{equation}
\end{lemma}

\begin{proof}
	First, we consider the two-dimensional case, i.e., $n = 2$, where the construction of $v_L$ is based on the rotated gradient operator, defined as $\nabla^{\perp} \coloneqq (-\partial_y, \partial_x)^\top$. Specifically, let $\eta \in H^1_0(T)$ denote the classical polynomial bubble function associated with $T$, obtained as the product of the corresponding hat functions and scaled such that $\eta(x_T) = 1$ at the barycenter $x_T$ of $T$. We then make the ansatz
	\begin{equation*}
		v_L = \nabla^{\perp} (\eta^2 p)
	\end{equation*}
	with the polynomial $p\in \mathbb{P}^{m - 1}(T)$  to be determined.
	This ansatz ensures that $\nabla \cdot v_L = 0$ on $T$ and $v_L = 0$ on $\partial T$. Using the ansatz and integrating by parts, condition~\cref{defbT} can be rewritten, for any $q \in (\mathbb{P}^m(T))^n$, as
	\begin{equation}
		\label{eq:condreformulated}
		\int_T \eta^2 p \nabla^{\perp} \cdot w \dx = -Lw.
	\end{equation}
	
	We decompose the polynomial space as $(\mathbb{P}^m(T))^n = \mathbb{G}^m(T) \oplus (\mathbb{G}^m(T))^{\perp}$, where $(\mathbb{G}^m(T))^{\perp}$ denotes the $L^2(T)$-orthogonal complement of $\mathbb{G}^m(T)$. Since for any $q \in \mathbb{G}^m(T)$ it holds that $\nabla^{\perp} \cdot q = Lq = 0$, it suffices to verify equation~\cref{eq:condreformulated} only for $q \in (\mathbb{G}^m(T))^{\perp}$. It was shown in \cite[Eq.~(2.7)]{daVeiga2015} that the operator $\nabla^{\perp} \cdot : (\mathbb{G}^m(T))^{\perp} \to \mathbb{P}^{m-1}(T)$ is an isomorphism. Denoting its inverse by $\mathcal{J} \colon \mathbb{P}^{m-1}(T) \to (\mathbb{G}^m(T))^{\perp}$, we can then equivalently rewrite equation~\cref{eq:condreformulated} as
	\begin{equation}
		\label{eq:eqforp}
		\int_T \eta^2 p r \dx = -L\mathcal{J}r,
	\end{equation}
	for all $r \in \mathbb{P}^{m - 1}(T)$.
	The left-hand side of this equation is a weighted $L^2(T)$-inner product and, therefore, the Riesz representation theorem ensures the existence of a unique $p$. Using the equivalence of norms in finite dimensions and a scaling argument, we obtain, after testing equation \cref{eq:eqforp} with $p$, that
	\begin{equation*}
		\| p \|_T \lesssim \| L \| \, \| 	\mathcal{J} \|_{L^2 \rightarrow L^2} \lesssim H \| L \|.
	\end{equation*}
	Using a classical inverse estimate for polynomials, the stability estimate \cref{estbT} for $v_L= \nabla^{\perp}(\eta^2_T p)$ can be concluded.
	
	In the three-dimensional case, i.e., $n = 3$, one proceeds similarly, but using the curl operator instead of the rotated gradient. Specifically, we make the ansatz $v_L = \nabla \times (\eta^2 p)$, where $\eta \in H^1_0(T)$ is the classical polynomial bubble function associated with $T$, scaled such that $\eta(x_T) = 1$ at the barycenter $x_T$, and $p$ is a polynomial in 
	$\mathbb{D}^{m-1}(T) \coloneqq \{ v \in (\mathbb{P}^{m-1}(T))^3 \with \div v = 0 \}$. 
	According to \cite[Eq.~(2.9)]{daVeiga2015}, the curl operator is an isomorphism from $(\mathbb{G}^m(T))^{\perp}$ to $\mathbb{D}^{m-1}(T)$, and the equation for $p \in \mathbb{D}^{m-1}(T)$ can be written for all $w \in (\mathbb{P}^m(T))^n$ as
	\begin{equation*}
		\int_T \eta^2_T p \cdot \nabla \times w \dx = L w,
	\end{equation*}
	which allows us to conclude the proof as in the two-dimensional case.
\end{proof}

\begin{proof}[Proof of \cref{lem:locladymod}]
	For the proof, we use an element-local version of Ladyzhenskaya's lemma. That is, for any $T \in \mathcal{T}_H$ and any $q \in L^2(T)$ with $\int_T q \, \mathrm{d}x = 0$, there exists $v \in (H^1_0(T))^n$ such that $\div v = q$ and 
	$\|\nabla v\|_T \lesssim \|q\|_T$, where the hidden constant depends only on the shape regularity of the mesh. This result can be derived, for example, from~\cite{bernardi16}, where the shape-dependence of the Ladyzhenskaya constant is investigated.  
	We further define a linear functional $L \colon (\mathbb{P}^m(T))^n \to \mathbb{R}$ by setting $L r \coloneqq \int_T v \cdot r \, \mathrm{d}x$ for all $r \in \mathbb{Q}^m(T)$, and $L(\mathbb{G}^m(T)) \coloneqq 0$. Using property~(3) of \cref{AssQm}, we obtain, for any polynomial $p \in (\mathbb{P}^m(T))^n$, with the decomposition $p=r+g$ into $r\in \mathbb{Q}^m(T)$ and $g \in \mathbb{G}^m(T)$, 
	\begin{align}
		\label{eq:boundopnorm}
		\| Lp \|_T = \| Lr \|_T \lesssim  \|p\|_T \, \|v\|_T \lesssim H \|p\|_T \, \|q\|_T,
	\end{align}
	where we used the Poincaré--Friedrichs inequality locally on $T$. Finally, we define $v_q \coloneqq v - v_L$, where $v_L \in (H^1_0(T))^n$ is given by \cref{LemCelWeights}. Thanks to~\cref{eq:boundopnorm}, which gives an estimate for $\|L\|$, cf.~\cref{estbT}, the function $v_q$ has the desired properties.
\end{proof}

\begin{proof}[Proof of \cref{LemBubbles}]
	Let us first construct the element bubble functions \(b_T\) satisfying \cref{eq:kdcells,eq:bubbleest2} for a given \(\mathbf{g}_T \in \mathbb{Q}^m(T)\). To this end, consider the linear functional \(L\colon (\mathbb{P}^m(T))^n \to \mathbb{R}\) defined by \(L(q) \coloneqq \int_T \mathbf{g}_T \cdot q \, dx\) for \(q \in \mathbb{Q}^m(T)\) and \(L(\mathbb{G}^m(T)) \coloneqq 0\).
	Using property (3) of \cref{AssQm}, as in the preceding proof, we obtain the bound \(\|L\| \lesssim \|\mathbf{g}_T\|_T\), where we recall the definition of the operator norm \(\|L\|\)  in \cref{estbT}. We then set \(b_T\) to be the function \(v_L \in (H^1_0(T))^n\) from \cref{LemCelWeights}, extended by \(0\) outside \(T\). This construction ensures that \(b_T\) satisfies~\cref{eq:bubbleest2,eq:kdcells}, and its divergence vanishes, so that \(b_T \in Z\).
	
	Next, we construct the face bubble function \(b_F\) satisfying \cref{eq:kdfaces,eq:bubbleest1} for a given \(g_F \in \mathbb{P}^m(F)\). Let \(b_F' \in H^1_0(\omega_F)\) denote the classical piecewise polynomial bubble function on \(F\), scaled so that \(b_F'(x_F) = 1\), where \(x_F\) is the barycenter of~\(F\). By the Riesz representation theorem, there exists \(\lambda \in \mathbb{P}^m(F)\) such that
	\begin{equation}
		\label{eq:condscalarbubble} \int_F b_F' \lambda \mu \ds = \int_F g_F \mu \ds
	\end{equation}
     for all $\mu \in \mathbb{P}^m (F)$.
   By scaling, we have \(\|b_F' \, \lambda\|_F \lesssim \|g_F\|_F\). Extending \(\lambda\) from \(F\) to \(\omega_F\) as a polynomial constant in the direction normal to \(F\), and multiplying the extension by \(b_F'\), we obtain a piecewise polynomial function \(v_F \in H^1_0(\omega_F)\) satisfying \(v_F = b_F' \, \lambda\) on \(F\). A scaling argument then shows that \(\|v_F\|_{\omega_F} + H \|\nabla v_F\|_{\omega_F} \lesssim H^{1/2} \, \|g_F\|_F\). Finally, the desired face bubble function \(b_F\) is defined locally on the two elements \(T_s\), \(s \in \{1,2\}\), that compose \(\omega_F\) as
	\begin{equation*}
		b_{F} |_{T_s} \coloneqq v_F |_{T_s}n_F  - b_{T_s} - v_{q,T_s} 
	\end{equation*}
    with \(n_F\) denoting the normal chosen on \(F\). Here, \(b_{T_s}\) are the element bubble functions constructed above, associated with \(T_s\) and corresponding to \(\mathbf{g}_{T_s} \in \mathbb{Q}^m(T_s)\) such that \(\int_{T_s} \mathbf{g}_{T_s} \cdot \mu \, dx = \int_{T_s} v_F n_F \cdot \mu \, dx\) for all \(\mu \in \mathbb{Q}^m(T_s)\). Moreover, \(v_{q,T_s} \in (H^1_0(T_s))^n\) are the local Ladyzhenskaya-type bubble functions on each \(T_s\), constructed in \cref{lem:locladymod} with \(q = \div(v_F n_F) - \frac{1}{|T_s|} \int_{T_s} \div(v_F n_F) \, dx\). It can be readily verified that this \(b_F\) satisfies the desired properties \cref{eq:bubbleest1,eq:kdfaces} and has piecewise constant divergence with respect to the mesh $\TH$, i.e., \(b_F \in Z\).
\end{proof}

We now turn to the fully discrete analogues of \cref{lem:locladymod,LemBubbles}. Note that the constant \(\gamma\) appearing in the lemma below is the same as in \cref{thm:fd}.

\begin{lemma}[Fully discrete Ladyzhenskaya-type bubble]\label{hlem:locladymod}
Assume that \(\tfrac{h}{H} \leq \gamma\) for some constant \(\gamma > 0\), depending only on the regularity of the meshes \(\mathcal{T}_H\) and~\(\mathcal{T}_h\). Then, for any \(T \in \mathcal{T}_H\) and any \(q_h \in Q_h(T)\) with \(\int_T q_h \, dx = 0\), there exists \(v_{q,h} \in V_h(T) \cap (H^1_0(T))^n\) satisfying 
\(b_T(v_{q,h}, \chi_h) = \int_T q_h \, \chi_h \dx\) for all \(\chi_h \in Q_h(T)\) and 
\(\int_T v_{q,h} \cdot p \dx = 0\) for all \(p \in \mathbb{Q}^m(T)\), 
and the following stability estimate holds:
    $$
		\|\nabla v_{q,h}\|_T \lesssim  \|q_h\|_T.$$
\end{lemma}

\begin{lemma}[Fully discrete bubble functions]\label{hLemBubbles}
Under the same assumption as in \cref{hlem:locladymod}, for any \(F \in \mathcal{F}_H^i\) and any \(g_F \in \mathbb{P}_m(F)\), there exists a fully discrete face bubble function \(b_{F,h} \in Z_h \cap (H^1_0(\omega_F))^n\) such that, for all $\bbmu \in M_H$,
	\begin{equation}
		\label{heq:kdfaces}
		c (b_{F, h}, \bbmu) = H\int_F g_F \mu \ds,
	\end{equation}
	and
\begin{equation*}
\|b_{F,h}\|_{\omega_F} \lesssim H^{1/2} \, \|g_F\|_F, \qquad 
\|\nabla b_{F,h}\|_{\omega_F} \lesssim H^{-1/2} \, \|g_F\|_F.
\end{equation*}

	Similarly, for any $T \in \mathcal{T}_H$ and any $\mathbf{g}_T \in   \mathbb Q^m (T) $, there exists a fully discrete element bubble function $b_{T, h} \in   Z_h \cap (H^1_0 (T))^n$ such that, for all $\bbmu \in M_H$,
	   \begin{equation}     \label{heq:kdcells} c (b_{T, h}, \bbmu) = \int_T \mathbf{g}_T \cdot     \bmu \dx,
        \end{equation}   
    and
    \begin{equation}
    \label{eq:stabbubbleTh}
        \|b_{T, h} \|_T  \lesssim  \| \mathbf{g}_T \|_T,\qquad \| \nabla b_{T, h} \|_T \lesssim H^{-1} \| \mathbf{g}_T \|_T
    \end{equation}
\end{lemma}

Before proving the  two lemmas above, we first establish a specific local inf--sup condition, from which they follow as corollaries.

\begin{lemma}[Local inf--sup result]
For any \(T \in \mathcal{T}_H\), there holds
\begin{equation} \label{infsupbhg}
\adjustlimits \inf_{q \in \mathbb Q^m(T)} \sup_{v_h \in Z_h^0(T)}
\frac{(v_h, q)_T}{(\|v_h\|_T^2 + H^2 \|\nabla v_h\|_T^2)^{1/2}\|q\|_T} \gtrsim 1,
\end{equation}
where \(Z_h^0 \coloneqq Z_h \cap (H^1_0(T))^n\).
\end{lemma}

\begin{proof}
For any \(q \in \mathbb Q^m(T)\), there exists, thanks to \cref{LemBubbles}, an element bubble function \(b_T \in H_0^1(T)^n\) satisfying \(\div b_T = 0\) on \(T\), \((b_T, q)_T = \|q\|_T^2\), and
\[
\|b_T\|_T + H \|\nabla b_T\|_T + H^2 |b_T|_{2,T} \lesssim \|q\|_T \,.
\]
The estimate for \(|b_T|_{2,T}\), although not stated explicitly in \cref{eq:bubbleest2}, follows directly from the construction of \(b_T\). Denoting \(V_h^0(T) \coloneqq V_h \cap (H_0^1(T))^n\) and \(X_h(T) = X_{H,h}|_T\), we construct a fine-scale discrete version of \(b_T\) by solving the following Stokes-like problem: seek \((v_h, r_h) \in V_h^0(T) \times X_h(T)\) such that
\begin{subequations}
    \begin{align}
 \qquad \qquad & (\nabla v_h, \nabla w_h)_T & - \quad & (r_h, \div w_h)_T & = \quad & (\nabla b_T, \nabla w_h)_T, & \quad \quad \label{eq:stokeslikesysa}\\
\qquad \qquad & (q_h, \div v_h)_T & & & = \quad & 0 & \quad \quad 
\end{align}
\end{subequations}
for all $(w_h,q_h ) \in V_h^0(T) \times X_h(T)$. Note that this problem is well-posed thanks to the assumed inf--sup condition \cref{infsuploc}.
Then, by standard approximation theory for Stokes problems, cf.~\cite{Girault1986}, we obtain the error estimate
\[ \| v_h - b_T \|_T \lesssim h^2 | b_T |_{2, T} \lesssim h^2H^{-2} \|q\|_T. \]   
If we now choose the constant \(\gamma\) so that the hidden constant in the previous inequality equals \(\frac{1}{2\gamma^2}\), then, using the assumption \(\frac{h}{H} \leq \gamma\), we obtain
\[
(v_h, q)_T \geq (b_T, q)_T - \|b_T - v_h\|_T \, \|q\|_T \geq \tfrac{1}{2} \, \|q\|_T^2 \,.
\]
Testing \cref{eq:stokeslikesysa} with \(v_h\) gives the stability estimate \(\|\nabla v_h\|_T \leq \|\nabla b_T\|_T\). Using this bound, we immediately obtain
\[
(\|v_h\|_T^2 + H^2 \|\nabla v_h\|_T^2)^{1/2} \lesssim H \|\nabla v_h\|_T \lesssim H \|\nabla b_T\|_T \lesssim \|q\|_T,
\]
and inf--sup condition \cref{infsupbhg} then follows noting that \(v_h \in Z_h^0(T)\).
\end{proof}
\begin{proof}[Proof of \cref{hlem:locladymod,hLemBubbles}]
Let us first construct the element bubble functions \(b_{T,h} \in Z_h^0(T)\) on any \(T \in \mathcal{T}_H\), as announced in \cref{hLemBubbles}. Given \(\mathbf{g}_T \in \mathbb Q^m(T)\), we define \((b_{T,h}, \blambda) \in Z_h^0(T) \times \mathbb Q^m(T)\) as the solution to
\begin{align*}
\qquad \qquad & (b_{T,h}, v_h)_T & + \quad & H^2 (\nabla b_{T,h}, \nabla v_h)_T & + \quad & (v_h, \blambda)_T & = \quad & 0, & \qquad \qquad \\
\qquad \qquad & (b_{T,h}, \bmu)_T & & & & &= \quad & (\mathbf{g}_T, \bmu)_T & \qquad \qquad
\end{align*}
for all $(v_h,\bmu) \in Z_h^0(T) \times \mathbb Q^m(T)$. Note that this saddle-point problem is well-posed thanks to the inf--sup condition \cref{infsupbhg}. Moreover, standard inf--sup theory, cf.~\cite[Cor.~4.2.1]{BoffiBrezziFortin2013}, implies the stability estimate \cref{eq:stabbubbleTh}. Extending \(b_{T,h}\) by \(0\) outside \(T\) then gives the desired fully discrete element bubble function.

Let us now construct \(v_{q,h} \in V_h^0(T)\) as announced in \cref{hlem:locladymod}. Given \(T \in \mathcal{T}_H\) and \(q_h \in Q_h(T)\), we define \((v_{q,h}', \xi_h) \in V_h^0(T) \times X_h(T)\) as the solution to
    \begin{subequations}
\begin{align}
\quad \qquad & (\nabla v_{q,h}', \nabla w_h)_T & - \quad & (\xi_h, \div w_h)_T & = \quad & 0, &  & \label{eq:stokesvh1} \\
\quad \qquad & (\chi_h, \div v_{q,h}')_T & & & = \quad & (\chi_h, q_h)_T &  & \label{eq:stokesvh2}
\end{align}
\end{subequations}
for all $(w_h, \chi_h) \in V_h^0(T) \times X_h(T)$.  This problem is well-posed thanks to the assumed inf--sup condition \cref{infsuploc}. Moreover, standard inf--sup theory, cf.~\cite[Cor.~4.2.1]{BoffiBrezziFortin2013}, implies the stability estimate \(\|\nabla  v_{q,h}'\|_T \lesssim \|q_h\|_T\). Let \(b_{T,h} \in Z_h^0(T)\) be the element bubble function constructed as above for \(\mathbf{g}_T \in \mathbb{Q}^m(T)\) such that \((\mathbf{g}_T, \bmu)_T = (v_{q,h}', \bmu)_T\) for all \(\bmu \in \mathbb{Q}^m(T)\). Then, we obtain
\[
\|\nabla b_{T,h}\|_T \lesssim H^{-1} \|\mathbf{g}_T\|_T \lesssim H^{-1} \| v_{q,h}'\|_T \lesssim \|\nabla v_{q,h}'\|_T,
\]
where we have used the stability estimate \cref{eq:stabbubbleTh} and the Poincaré–Friedrichs inequality on \(T\). The function defined by \(v_{q,h} \coloneqq v_{q,h}' - b_{T,h}\) then satisfies all the properties stated in \cref{hlem:locladymod}.

Finally, the fully discrete face bubble function \(b_{F,h}\) from \cref{hLemBubbles} can be constructed analogously to the continuous case in Lemma \ref{LemBubbles}, using the previously defined fully discrete element bubbles \(b_{T,h}\) and the fully discrete functions \(v_{q,h}\) supported on the two mesh elements adjacent to \(F\).
\end{proof}

\section{Inf--sup stability of modified Taylor--Hood}\label{sec:InfSupTHmod}

\begin{lemma}[Modified Taylor--Hood]\label{InfSupTHmod}
The finite element pair \(V_h / Q_h\), with \(V_h=\{v_h\in (H^1_0(\Omega))^n:v_h|_K\in\mathbb{P}^{k+1}(K)\ \forall K\in\Th\}\) and \(Q_h=(\oplus_{T\in\TH}Q_h(T))\cap L^2_0(\Omega)\) with \(Q_h(T)=\{q_h\in H^1(T):q_h|_K\in\mathbb{P}^{k}(T)\ \forall K\in\Th(T)\}\)  on each coarse element \(T \in \mathcal{T}_H\), satisfies the following inf--sup property:
	\[ \adjustlimits \inf_{p_h \in Q_h} \sup_{v_h \in V_h} \frac{b(v_h,p_h )}{ \|\nabla  v_h \|_{\Omega}\| p_h \|_{\Omega}} \gtrsim 1.\]
\end{lemma}

\begin{proof}
Take any $p_h \in Q_h$ and decompose it as $p_h = p_h^0 + \bar{p}_h$ with $p_h^0 \in X_{H,h}$ and $\bar{p}_h \in \mathbb{P}^0(\TH)$. By the inf--sup stability of the pair $V_h(T)/Q_h(T)$ on each $T \in \TH$, and recalling that $\int_T p_h^0 \, \mathrm{d}x = 0$, there exists $v_h^T \in V_h(T)$ such that
	\[ (p_h^0, v_h^T)_T = \| p_h^0 \|_T^2,\qquad \|\nabla  v_h^T
	\|_{T} \lesssim \| p_h^0 \|_T  \]
with the hidden constant related to the inf--sup constant of this finite element pair, which depends only on the regularity of the fine and coarse meshes (through the shape of $T$). We define $v_h^0 \in V_h$ such that $v_h^0|_T = v_h^T$ for every $T \in \TH$.

	To handle $\bar{p}_h$, we first recall that there exists $v \in (H^1_0(\Omega))^n$ such that $\div v = \bar{p}_h$ and $\|\nabla v\|_{\Omega} \lesssim \|\bar{p}_h\|_\Omega$. Let $\bar{v}_h \in V_h$ be a suitably chosen interpolant of $v$ (to be constructed) such that $\int_E \bar{v}_h\ds  = \int_E v\ds $ on every interior facet $E$ of the coarse mesh $\TH$, and $\|\nabla \bar{v}_h\|_{\Omega} \lesssim \|\nabla v\|_{\Omega} \lesssim \|p_h^0\|_\Omega$. Then, on every element $T \in \TH$,
	\[ b_T (\bar{p}_h, \bar{v}_h) = \bar{p}_h |_T  \int_{\partial
		T} \bar{v}_h \cdot n = \bar{p}_h |_T  \int_{\partial T} v \cdot
	n = \| \bar{p}_h \|_T^2.\]
	Let $v_h = v_h^0 + \lambda \bar{v}_h$, where $\lambda = \operatorname{sgn}\big(b(p_h^0, \bar{v}_h)\big)$, with the sign function $\operatorname{sgn}(\cdot)$. Then,
	\[ b(p_h, v_h) = b(p_h^0 + \bar{p}_h,	v_h^0 + \lambda \bar{v}_h) = \| p_h^0 \|_{\Omega}^2 + \lambda \|	\bar{p}_h \|_{\Omega}^2 + \lambda  b(p_h^0, 
	\bar{v}_h)  \geq \| p_h\|_{\Omega}^2.  \]
	We also have
	\[ | v_h |_{1, \Omega} \leq | v_h^0 |_{1,\Omega} + | \bar{v}_h |_{1,\Omega}  \lesssim \| p_h^0 \|_{\Omega} + \| \bar{p}_h \|_{\Omega}  \lesssim  \| p_h \|_{\Omega}, \]
	which implies the announced inf--sup inequality.
	
	It remains to construct an interpolant $\bar{v}_h \in V_h$ for $v \in H^1_0(\Omega)^d$, which preserves the integrals on the faces of the coarse mesh and satisfies $|\bar{v}_h|_{1,\Omega} \lesssim |v|_{1,\Omega}$. We describe this construction in the two-dimensional case (i.e., $n = 2$). We start from the usual Cl\'ement-type interpolant $I_h v \in V_h$, which satisfies $\|\nabla I_h v\|_{K} \lesssim \|\nabla v\|_{\omega_K}$ and $\|v - I_h v\|_E \lesssim \sqrt{h}\,\|\nabla v\|_{\omega_E}$ on any fine mesh element $K \in \Th$ and any interior edge $E$ of the fine mesh $\Th$. Let $\mathcal{E}_{H,h}$ denote the set of edges of $\Th$ lying on an edge of the coarse mesh $\TH$. For any $E \in \mathcal{E}_{H,h}$, let $\phi_E \in V_h$ denote the piecewise quadratic polynomial on $\Th$ equal to 1 at the midpoint of $E$ and vanishing at all other edge midpoints and nodes of the fine mesh. We set
	\[ \bar{v}_h = I_hv + \sum_{E \in \mathcal{E}_{H, h}} \frac{\int_E (v - I_hv)\ds} {\int_E \phi_E\ds} \phi_E \]
	so that $\int_E \bar{v}_h \ds = \int_E v \ds$ for all $E \in \mathcal{E}_{H,h}$. Note that $\bar{v}_h \in V_h$, since the finite element space is of at least quadratic order. Moreover, we also have
\[
\|\nabla \bar{v}_h\|_{\Omega}^2 \lesssim \|\nabla I_h v\|_{\Omega}^2 + \sum_{E \in \mathcal{E}_{H,h}} \hspace{
-1ex}\|\nabla v\|_{\omega_E}^2 \, \|\nabla \phi_E\|_{\omega_E}^2
\lesssim \|\nabla v\|_{\Omega}^2 + \sum_{E \in \mathcal{E}_{H,h}} \hspace{
-1ex}\|\nabla v\|_{\omega_E}^2
\lesssim \|\nabla v\|_{\Omega}^2.
\]
	Thus, the interpolant $\bar{v}_h$ has indeed the announced properties.
	
	The proof in the three-dimensional case (i.e., $n = 3$) is similar for polynomial degrees $k \geq 2$. In that case, it suffices to replace the quadratic functions $\phi_e$ by cubic functions associated with the barycenters of faces in $\mathcal{E}_{H,h}$. The construction for $n = 3$ and $k = 1$ is more involved, since the cubic functions are no longer contained in $V_h$. Nevertheless, the desired properties can still be achieved using suitable combinations of quadratic basis functions associated with edges rather than faces. We omit the details of this more tedious construction.
\end{proof}
\bibliographystyle{alpha}
\bibliography{bib}

\end{document}